\input amstex
\documentstyle{amsppt}
\magnification=\magstep1 \catcode`\@=11
\catcode`\@=\active

\voffset-.5in

\document

\NoRunningHeads
\def\C{\Bbb C}

\def\H{\Bbb H}

\def\R{\Bbb R}
\def\Z{\Bbb Z}
\def\aa{\frak a}

\def\TTT{\Cal T}
\def\PPP{\Cal P}

\def\wh{\widehat}
\def\wt{\widetilde}
\def\wmu{\widetilde\mu}
\def\vf{\varphi}
\def\ve{\varepsilon}

\def\o{\overline}
\def\tphi{\;^t\phi}

\def\e{\epsilon}
\def\k{\kappa}
\def\pa{\parallel}

\def\bG{{\bold G}}
\def\bGL{{\bold {GL}}}
\def\bT{{\bold T}}
\def\bP{{\bold P}}
\def\bM{{\bold M}}
\def\bN{{\bold N}}
\def\bA{{\bold A}}
\def\bB{{\bold B}}
\def\bS{{\bold S}}
\def\b1{{\bold 1}}

\heading Langlands Classification for L-parameters\endheading

\centerline{by Allan J. Silberger and Ernst-Wilhelm Zink}
\medskip
\noindent {\bf Abstract:} {\it Let $F$ be a non-archimedean local field and $G=\bG(F)$ the group of $F$-rational points
of a connected reductive $F$-group. Then we have the Langlands classification of complex irreducible admissible
representations $\pi$ of $G$ in terms of triples $(P,\sigma,\nu)$ where $P\subset G$ is a standard $F$-parabolic
subgroup, $\sigma$ is an irreducible tempered representation of the standard Levi-group $M_P$ and
$\nu\in \R\otimes X^*(M_P)$ is regular with respect to $P.$ Now we consider Langlands' L-parameters $[\phi]$
which conjecturally will serve as a system of parameters for the representations $\pi$ and which are
(roughly speaking) equivalence classes of representations $\phi$ of the absolute Galois group
$\Gamma=\text{Gal}(\overline{F}|F)$ with image in Langlands' L-group $\,^LG$, and we classify the possible $[\phi]$
in terms of triples $(P,[\tphi],\nu)$ where the data $(P,\nu)$ are the same as in the Langlands classification of representations and where $[\tphi]$ is a tempered L-parameter of $M_P.$}
\bigskip
In analogy with an improved version of A.~Silberger's 1978 paper {\bf [S]}: on the Langlands Quotient Theorem for p-adic Groups, which reduces the classification of irreducible representations to the classification of the tempered representations, we present a similar reduction in terms of Langlands' L-parameters which are conjecturally in a 1-1 correspondence with the L-packets of irreducible representations. Depending on that we add some remarks on {\bf refined} L-parameters, where an additional parameter for the constituents of a single L-packet has to be added.

\noindent Let $G=\bG(F)$ be a connected reductive group over a p-adic field and let $\Gamma=\text{Gal}(\overline F|F)$ be the absolute Galois group, where $\overline F$ denotes a separable algebraic closure of $F.$ An L-parameter $[\phi]$ for $G$ is roughly speaking an equivalence class of L-homorphisms $\phi$ which represent $\Gamma$ in Langlands' L-group $\,^LG$ (for details see section {\bf 3}).
Just as the Langlands classification of irreducible representations $\pi$ of $G$ is performed in terms of standard triples $(P,\sigma,\nu)\leftrightarrow\pi,$ where $\sigma$ is a tempered representation of the standard $F$-Levi subgroup $M_P,$ (see the reminder in section 1), {\it we want to present now a Langlands classification of L-parameters $[\phi]$ of $G$ in terms of standard triples $(P,[\tphi],\nu)$ where $(P,\nu)$ are the same data as in the case of representations and where $[\tphi]$ is a tempered L-parameter of $M_P$} (see 3.5 and 4.6). The local Langlands conjecture expects the set $Irr(G)$ of equivalence classes of irreducible representations as the disjoint union of L-packets  of finite size:
$$   Irr(G) = \bigcup_{[\phi]\in\Phi(G)} \Pi_{[\phi]},$$
where $\Phi(G)$ denotes the set of all L-parameters. The Langlands classification of L-parameters suggests
that all representations $\pi\in \Pi_{[\phi]}$ in a fixed packet must have the same data $(P,\nu)$, because these data are already determined by the L-parameter $[\phi]$, whereas $\sigma$ may vary over the tempered L-packet $\Pi_{[\tphi]}\subset Irr(M_P).$ (For this and other consequences see section {\bf 7}).

\noindent The Langlands classification of L-parameters has been sketched before in papers of
J.Arthur [A] and implicitly in V.Heiermann [H] which have been stimulating for our approach. (In [A] it is only sketched as a first reduction to the main objective of that paper, to study the tempered L-packets). Anyway it seemed to us that these papers passed over several details which are worthwhile to be spelled out as we try to do in the following remarks. The contents of this paper is:

{\bf 1.} Langlands classification of irreducible representations - a reminder.

{\bf 2.} Langlands' L-group, a review with additional remarks. The main objective is to introduce relevant semistandard parabolic and Levi subgroups in $\,^LG$ which are in a natural 1-1 correspondence with semistandard $F$-parabolic and $F$-Levi subgroups in $G=\bG(F)$ (see 2.5.4, 2.5.5) and to identify the relative Weyl-group $\,_FW$ in terms of the dual group $\wh G$ (see 2.4.2 and 2.5.4(ii)) to make this 1-1 correspondence equivariant under the Weyl-group actions. It is then a main tool in section 6. A reader who is willing to accept this can skip the section.

{\bf 3.} L-parameters and their Langlands classification: Here we repeat the basic definitions on L-homomorphisms and recall the local Langlands conjecture. Related to that the notion of a tempered L-parameter comes up, and - having in mind 1.4 - this leads to a first, still preliminary, statement of what a Langlands classification of L-parameters should do (see 3.5).

{\bf 4.} Twisting L-homomorphisms and a precise statement of 3.5: The Langlands classification of irreducible representations uses the possibility of twisting a representation of a Levi subgroup $M=\bM(F)$ by a positive real valued unramified character of $M.$ These characters are in natural correspondence with the elements of a real vector space $\aa_M^*$ (see 1.2(2)). Now we interpret $\aa_M^*$ as the hyperbolic elements in the central torus $Z(\,^LM)^0$  of the L-group (see 4.3), and combine this with the possibility of twisting an L-homomorphism (which takes its values in $\,^LM$) by an element of $Z(\,^LM)^0$ (see 4.5). This leads us to 4.6 as a precise formulation of 3.5. To motivate this we remark in 4.7 the conjectural behaviour of the Langlands map with respect to unramified twists.

{\bf 5.} Two basic invariants associated to an L-homomorphism: In 5.3 we follow {\bf [H],5.1} to obtain a
map $\phi\mapsto z(\phi)$ which associates with the L-homomorphism $\phi$ a hyperbolic element in the central torus of the reductive group $C_{\wh G}(Im\,\phi),$ and in 5.5 we determine in $\;^LG$ a distinguished Levi-subgroup $L(\;^LG)_\phi$ depending on $\phi.$ The natural map $\phi\mapsto \left(z(\phi),\;L(\;^LG)_\phi\right)$ commutes with $\wh G$-conjugation.

{\bf 6.} On the proof of 3.5 / 4.6. It is obvious how to obtain an L-parameter out of an L-parameter standard triple. Now we use {\bf 5.} and the information obtained in {\bf 2.} on relevant groups, in order to construct the converse map which, beginning from an L-parameter $[\phi],$ determines the corresponding L-parameter standard triple.

{\bf 7.} A remark on refined L-parameters. As has been mentioned above we can expect a  natural identification $\Pi_{[\phi]}\leftrightarrow \Pi_{[\,^t\phi]_{M_P}}$ between a general L-packet and a tempered L-packet if the L-parameter $[\phi]$ is given by the L-parameter standard triple $(P,[\tphi],\nu).$ (The possibility of such an identification has been predicted in {\bf [A],p.201}.) In particular the constituents of the packet $\Pi_{[\phi]}$ allow the same parameters as those of the tempered packet
$\Pi_{[\,^t\phi]_M}.$ They are (conjecturally) given in terms of a component group $\overline S_{\,^t\phi}(M)$ which, as we will prove in {\bf 7.1}, is equal to the component group $\overline S_\phi(G)$
which carries the parameters for the packet $\Pi_{[\phi]}.$ Here we thank Ahmed Moussaoui for referring us to Lemma 1.1 of [A 99].

{\bf 8.} As examples we recover the Langlands classification of L-parameters for $GL_n$ and for its inner forms.

\noindent Our base field $F$ will be a locally compact non-archimedean field, and we will denote

$q=\#\kappa_F$ the order of the finite residue field of $F$,

$p$ the characteristic of $\kappa_F,$

$\nu_F:F \twoheadrightarrow \Z\cup\{\infty\}$ the exponent on $F$,

$|a|_F= q^{-\nu_F(a)}$ the normalized absolute value on $F.$

\noindent For short we will speak of $F$ as of a p-adic field.

With respect to algebraic groups we use the following facts:

\noindent The $F$-group $\bG$ is completely determined (up to $F$-isomorphism) by the group $\bG(\overline{F})$ of $\overline F$-rational points, where $\overline F|F$  is a separable closure. Therefore we will identify $\bG=\bG(\overline F).$  {\bf [Sp1], 3.1}

\noindent  A closed subvariety $Z$ of an $F$-variety $V$ is by itself an $F$-variety if and only if $Z$ is defined over $\overline F$ and is $\Gamma= \text{Gal}(\overline F|F)$-stable.  {\bf [Bo2], AG 14.4}

\noindent  If $\bG|F$ is a connected reductive group, then it is the inner form of a quasisplit group.  {\bf [Sp1], 3.2, [Sp2],16.4.9.}


We thank Anne-Marie Aubert for reading the manuscript and for her helpful remarks and suggestions. And the second named author wants to thank Peter Schneider for a useful discussion concerning 3.2 .
\bigskip

{\bf 1. Langlands classification of irreducible representations - a reminder}
\bigskip
{\bf 1.1 Parabolic and Levi subgroups}

Let $\bG|F$ be a connected reductive group which is defined over a
p-adic field $F.$
In $\bG$ we fix a maximal $F$-split torus $\bA_0$ and
denote $\bM_0$ its centralizer which is a minimal $F$-Levi subgroup in
$\bG.$ Moreover we fix a minimal $F$-parabolic subgroup $\bP_0$ which
has $\bM_0$ as one of its Levi subgroups, hence $\bP_0=\bM_0\bN_0$ where
$\bN_0$ is the unipotent radical of $\bP_0$. Then the {\bf semi-standard} and {\bf standard}
$F$-parabolic subgroups $\bP$ of $\bG$ are those which contain $\bM_0$ or even $\bP_0$ resp. In
a semi-standard group $\bP$ we have uniquely determined subgroups
$$\bN_P,\quad \bM_P,\quad \bA_P,$$
which are the unipotent radical, the {\bf uniquely determined} Levi subgroup which contains
$\bM_0$ and the maximal $F$-split torus in the center of $\bM_P,$ resp.
All algebraic groups $\bG$ we are going to work with
are defined over $F$ and we use ordinary letters $G=\bG(F)$ to denote the
corresponding groups of $F$-rational points.
The surjective map
$$   P\mapsto  M_P$$
from the set of semi-standard $F$-parabolic subgroups onto the set of semi-standard $F$-Levi subgroups
in $G$ is of particular importance for us.
\bigskip
{\bf 1.2 Unramified characters}

If $M=\bM(F)$ is a p-adic group, which in our context always will denote a Levi-subgroup
(of a parabolic subgroup) of $G$, then an {\bf unramified character} $\chi:M\rightarrow \C^\times$
is by definition a continuous homomorphism which is trivial on all compact subgroups of $M.$
We denote $X^*(M)=X^*(\bM)_F$ the
group of rational characters $\chi:\bM\rightarrow {\bold {GL}}_1$ which are defined over $F$
and therefore give rise to $\chi:M\rightarrow F^\times.$ By convention the abelian groups
$X^*(M)$ are written additively and therefore we sometimes write $\chi\in X^*(M)$ exponentially.
The group $X_{ur}(M)$ of unramified characters can be given in terms of $F$-rational characters
as follows:
$$ \C\otimes X^*(M) \twoheadrightarrow \C^\times\otimes_\Z X^*(M) \twoheadrightarrow X_{ur}(M), \tag{1}$$
where the first arrow takes $s\otimes\chi$ to $(q^s)^{\otimes \chi},$  $q=\#$residue field of $F$,

\noindent and the second arrow takes
$\lambda^{\otimes\chi}$ to the unramified character $\{m\mapsto \lambda^{-\nu_F(\chi(m))}\}$.
{\bf Note that the second surjection has finite kernel} because we have here two complex tori of the same rank.

Using the normalized absolute value on $F$ the combined map is:
$$  s\otimes\chi\mapsto |\chi|_F^s:=\{m\mapsto |\chi(m)|_F^s\}.$$
Note that $|\chi|_F^s$ is unramified because the rational character $\chi:M\rightarrow F^\times$ takes
compact subgroups of $M$ into the subgroup
of units in $F.$

\noindent The subgroup $X_{ur}(M)_+$ of positive real valued unramified characters
$\alpha:M\rightarrow \R_+^\times$
comes via the subgroups
$$  \aa_M^*:=\R\otimes X^*(M)\cong (\R_+)^\times\otimes X^*(M) \cong X_{ur}(M)_+ \tag{2}$$
(here we use that $s\in\R \mapsto q^s\in\R_+$ is a bijection) which we will write as $$\nu\in\aa_M^*
\mapsto \chi_\nu.$$ Contrary to (1) the restricted map $s\in\R \mapsto q^s$ turns
(2) into an isomorphism. If $M$ is semi-standard, i.e. $M\supseteq M_0$ then via restriction from $M$ to $M_0$ we have $\aa_M^*\subseteq \aa_{M_0}^*$
and we consider
them as {\bf euclidean spaces} with respect to a {\bf fixed Weyl-group invariant pairing $<\;,\;>$ on $\aa_{M_0}^*$}. (Here of course we think of the relative Weyl group $\;_FW= W(A_0,G),$ contrary to the absolute Weyl group $W(\bT,\bG)$ which will occur in section 2.)
\bigskip
{\bf 1.3 A standard triple} $(P,\sigma,\nu)$

\noindent for the reductive p-adic group $G=\bG(F)$ consists of:

{\bf a)} $P=\bP(F)$ a standard $F$-parabolic subgroup,

{\bf b)} an irreducible {\bf tempered} representation of the standard Levi subgroup $M_P,$

{\bf c)} a real parameter $\nu\in \aa_P^*:=\aa_{M_P}^*$ which is {\bf regular,} that is properly contained in the chamber which is determined by $P$. More precisely:
Let $\Delta(P)\subset X^*(A_P)$ be the simple roots for the adjoint
action of the torus $A_P$ on $Lie(N_P)$. Via restriction we may
identify $\aa_P^*=\R\otimes X^*(M_P)=\R\otimes X^*(A_P)\supset
\Delta(P),$ and the condition for $\nu$ to be regular is:
$$   \langle \nu,\alpha\rangle \; > 0\quad \forall \;\alpha\in \Delta(P).$$

{\bf d)} To put things into perspective we still mention that

$\aa_M^*=\Delta(A_0,M\cap P_0)^\perp\hookrightarrow \aa_{M_0}^*,$ and

$ \text{res}^{A_0}_{A_M}:\aa_{M_0}^*\rightarrow\aa_M^*\hookrightarrow\aa_{M_0}^*$ is an orthogonal projection which induces a bijection $\Delta(P_0)-\Delta(A_0,M\cap P_0)\leftrightarrow \Delta(P)\subset \aa_M^*\cap{\aa_G^*}^\perp.$

\noindent So, if $\nu\in\aa_{M'}^*=\Delta(A_0,M'\cap P_0)^\perp$ for some standard $M'\supset M$, then $\Delta(A_0,M'\cap P_0)-\Delta(A_0,M\cap P_0)$ gives rise to elements $\alpha\in \Delta(P)$ such that $\langle \nu,\alpha\rangle =0,$ hence $\nu$ cannot be regular. Therefore $\nu\in\aa_P^*$ regular implies that $\nu\notin \aa_{M'}^*$ for larger standard Levi groups $M'\supset M_P$ which means that $\chi_\nu\in X_{ur}(M_P)_+$ cannot be extended to $M'.$
\bigskip
\proclaim {1.4 Theorem (Langlands classification)} In $G$ we fix the data $(P_0,M_0)$ and the euclidean pairing $<\;,\;>$ on $\aa_{M_0}^*.$ Then there exists a natural bijection
$$  \{(P,\sigma,\nu)\} \longleftrightarrow Irr(G)\tag{3}$$
between the set of all standard triples and the set of (equivalence classes of) irreducible representations of $G.$ It takes a triple $(P,\sigma,\nu)$ to the representation
$$  \pi=j(P,\sigma,\nu):= j(i_{G,P}(\sigma\otimes\chi_\nu)),$$
where $i_{G,P}$ denotes the normalized parabolic induction from $M_P$ to $G$, and $j$ denotes the uniquely determined irreducible quotient of that induction.\endproclaim
\proclaim {1.4* A second realization}

\noindent A different bijection between standard triples and irreducible representations can be realized as follows:
$$  \pi=s(P,\sigma,\nu):= s(i_{G,P}(\sigma\otimes\chi_{-\nu})),\tag{3*}$$
where $s$ denotes the uniquely determined irreducible subrepresentation of that induction.

\noindent The two realizations are connected via forming the contragredient representation:
$$  s(P,\sigma,\nu) = \wt j(P,\wt\sigma,\nu),\qquad j(P,\sigma,\nu)= \wt s(P,\wt\sigma,\nu).$$
\endproclaim

\bigskip
{\bf 2. Langlands' L-group, a review}

\bigskip
{\bf 2.1 Defining the L-group}

Let $\bG|\overline F$ be a connected reductive group over a separable closure, and fix
$$  \bG\supset\bB\supset\bT $$
a Borel subgroup and a maximal torus. This determines a {\bf based root datum}
$$ \psi_0(\bG)=\{ X^*(\bT),\quad \Delta^*(\bT,\bB),\quad X_*(\bT),\quad \Delta_*(\bT,\bB)\},\tag{1}$$
consisting of the lattice of all rational characters of $\bT$, the set of simple roots with respect to $\bB$, and dually the lattice of cocharacters and the set of simple coroots.
\bigskip
Then there is a well defined {\bf complex dual group} $(\wh G,\wh B,\wh T)|\C,$ such that
$$ \psi_0(\wh G)=\{X_*(\wh T),\quad \Delta_*(\wh T,\wh B),\quad X^*(\wh T),\quad \Delta^*(\wh T,\wh B)\}\tag{2}$$
identify with the object placed at the same position of the upper row (1).

Now assume that $\bG|F$ is defined over a p-adic field. Then the Galois
group $\Gamma:=Gal(\overline F|F)$ acts on the objects (1), but in a modified way if $  \bG\supset\bB\supset\bT $ is not quasisplit. Nevertheless choosing a splitting, i.e. a set of generators $x_\alpha\in \bG_\alpha$ for the root subgroups corresponding to the simple roots $\alpha\in\Delta:= \Delta^*(\bT,\bB),$  we obtain a homomorphism ({\bf [Bo1], 1.2, 1.3} and {\bf [Sp1], 3.2, 2.13})
$$  \mu=\mu_\bG:\Gamma\rightarrow Aut\;\psi_0(\bG) \cong Aut(\bG,\bB,\bT,\{x_\alpha\}_{\alpha\in\Delta})\subset Aut(\bG,\bB,\bT),\tag{3}$$
and we speak here of the $\mu(\Gamma)$-action which has to be distinguished from the usual $\Gamma$-action on $\bG=\bG(\overline F).$ We recall the modified Galois action on $\psi_0(\bG):$
$$ \gamma\cdot_\mu\chi= \gamma\chi\gamma^{-1}\circ Int(a_\gamma)\in X^*(\bT),\tag{3a}$$
$$ \gamma\cdot_\mu\chi^\vee = Int(a_\gamma^{-1})\circ \gamma\chi^\vee\gamma^{-1}\in X_*(\bT),\tag{3b}$$
for $a_\gamma\in\bG$ such that $a_\gamma(\bB,\bT)a_\gamma^{-1} = (\gamma(\bB),\gamma(\bT)),$ hence the coset $a_\gamma\bT$ is uniquely determined.

\noindent {\bf Note} that actually $\mu$ depends on the triple $(\bG,\bB,\bT)$. And $\mu(\gamma)\in Aut\;\psi_0(\bG)$ lifts to $\wt \mu(\gamma)\in Aut(\bG,\bB,\bT)$ such that
$$ \psi_0(\wt\mu(\gamma))^* =\mu(\gamma^{-1})\;\in Aut\;X^*(\bT),\qquad \psi_0(\wt\mu(\gamma))_*= \mu(\gamma)\;\in Aut\;X_*(\bT).$$
For instance on $\bT =\overline F^\times\otimes X_*(\bT)$ the lifted $\mu(\Gamma)$-action is only on the second factor. The natural pairing
$$  X^*(\bT)\times X_*(\bT) \rightarrow Hom(GL_1,GL_1)=\Z,\qquad \langle\gamma\cdot_\mu\chi,\; \gamma\cdot_\mu\chi^\vee\rangle =\langle\chi,\; \chi^\vee\rangle,\tag{3c}$$
which is obtained by applying characters to cocharacters, is $\mu(\Gamma)$-invariant.

\noindent By transport of structure, $\Gamma$ acts on (2), $\mu_{\wh G}:\Gamma\rightarrow Aut(\psi_0(\wh G)),$
and this lifts to  an action $\wt\mu_{\wh G}(\Gamma)$ on $(\wh G,\wh B,\wh T),$ where again we have to fix a splitting $\wh x_{\wh\alpha}\in \wh G_{\wh\alpha}$ for all $\wh\alpha\in\Delta^*(\wh T,\wh B)$ in order to make this action unique. Otherwise it is only determined up to conjugation by elements from $\wh T.$ {\it Here we simply speak of the $\Gamma$-action on $\wh G$ and omit the notation $\wt\mu_{\wh G}.$ }

\noindent{\bf Langlands' L-group} is the semidirect product
$$  \;^LG:= \wh G\rtimes W_F$$
where the Weil group $W_F\subset \Gamma$ acts on the dual group via $\wt\mu_{\wh G}:\Gamma\rightarrow Aut(\wh G,\wh B,\wh T).$

\bigskip

{\bf 2.2 Remarks on the Weyl groups}

\noindent Using the natural bijection between roots and coroots we obtain the identification
$$\Delta^*(\bT,\bB)\leftrightarrow  \Delta_*(\bT,\bB)= \Delta^*(\wh T,\wh B)$$
between {\bf simple roots of the original group} and {\bf simple roots of the dual group} which we will write $\alpha\mapsto\wh\alpha.$ The Weyl groups
$$ W:=W(\bT,\bG)\subset Aut(X^*(\bT)),\qquad \wh W:=W(\wh T,\wh G)\subset Aut(X^*(\wh T))$$
are generated by the reflections $s_\alpha:$ $s_\alpha(x)=x-<x,\wh \alpha>\alpha$ and $s_{\wh \alpha}$  resp., where the pairing
$$ X^*(\bT)\times X^*(\wh T)\rightarrow\Z$$
is defined via $X^*(\wh T)=X_*(\bT).$ From $<\alpha,\wh\alpha>=2$ we have $\langle s_\alpha(x), s_{\wh \alpha}(y)\rangle =\langle x,y\rangle$ for all $x\in X^*(\bT),$ $y\in X^*(\wh T),$ and therefore an isomorphism
$$   W(\bT,\bG)\cong W(\wh T,\wh G),\qquad w\mapsto \wh w:=\,^tw^{-1},\tag{4}$$
such that $\wh{s_\alpha}=s_{\wh\alpha} \in W(\wh T,\wh G)$.

\noindent Now the $\mu(\Gamma)$-action on $X^*(\bT)$ and $X_*(\bT)= X^*(\wh T)$ together with the property {\bf (3c)}, induces a {\bf Galois-action on the Weyl groups} (via permutations of the generating reflections):
$$  w\mapsto \mu(\gamma)\cdot w\cdot \mu(\gamma)^{-1}\in Aut(X^*(\bT))\qquad  \wh w\mapsto \mu(\gamma)\cdot \wh w\cdot \mu(\gamma)^{-1}\in Aut(X_*(\bT)),$$
such that, via transport of the Galois structure, the isomorphism {\bf (4)}
becomes $\mu(\Gamma)$-$\mu_{\wh G}(\Gamma)$-equivariant.

\proclaim{2.2.1 Lemma} If $\bG,$ $\bG'$ are $F$-groups which are inner forms of each other and if
$$  f:(\bG,\bB,\bT)@>\sim >> (\bG',\bB',\bT')$$
is a corresponding $\overline F$-isomorphism, then the induced isomorphisms $f^*:X^*(\bT')@> \sim >>X^*(\bT)$ and
$$   w\in W(\bT,\bG) \mapsto (f^*)^{-1}\cdot w\cdot f^* \;\in W(\bT',\bG')$$
will be $\mu(\Gamma)$-equivariant, which means that $\mu_\bG(\gamma)\cdot w\cdot \mu_\bG(\gamma)^{-1}$ will be taken to
$\mu_{\bG'}(\gamma)\cdot [(f^*)^{-1}\cdot w\cdot f^*]\cdot \mu_{\bG'}(\gamma)^{-1}.$
\endproclaim
\demo{Proof} Because $f$ takes $(\bB,\bT)$ to $(\bB',\bT'),$ the map $f^*:X^*(\bT')\rightarrow X^*(\bT)$ will take simple roots $\alpha'\in \Delta(\bT',\bB')$ into simple roots $\alpha= f^*(\alpha')\in \Delta(\bT,\bB).$ Similar the map
$$  x\in Aut(X^*(\bT))\mapsto (f^*)^{-1}\cdot x\cdot f^* \in Aut(X^*(\bT'))$$
will map the reflection $s_\alpha$ onto $s_{\alpha'},$ and therefore it maps $W(\bT,\bG)$ onto $W(\bT',\bG').$ Moreover since our groups are inner forms of each other we know from {\bf [Bo1],1.3} that
$$   f^*\cdot \mu_{\bG'}(\gamma) = \mu_\bG(\gamma)\cdot f^*:\quad X^*(\bT')\rightarrow X^*(\bT),$$
for all $\gamma\in \Gamma,$ and therefore the result follows.
\qed\enddemo
\definition{Remark} In particular we see that the $\mu_{\bG'}(\Gamma)$-action on $X^*(\bT')$ and on $W(\bT',\bG')$ will be trivial if $\bG'|F$ is an inner form of a split group. This of course is not the $\Gamma$-action.
\enddefinition

\noindent Furthermore we need to switch to the description $W(\bT,\bG)\cong N_\bG(\bT)/\bT,\; w\mapsto n_w$ such that:
$$   w(\chi)= n_w\chi n_w^{-1},\quad\forall \chi\in X^*(\bT).$$
Here we prove:
\proclaim{2.2.2 Lemma} If $w = Int(n_w)\in Aut(X^*(\bT))$ and if $\wt\mu(\gamma)\in Aut(\bG,\bB,\bT)$ is a lift of $\mu(\gamma)=\mu_{\bG}(\gamma)\in Aut(\psi_0(\bG)),$ then:

{\bf (i)} $\mu(\gamma)\cdot w\cdot\mu(\gamma)^{-1}= Int(\wt\mu(\gamma)(n_w)),$

{\bf (ii)} $\mu(\gamma)\cdot w\cdot\mu(\gamma)^{-1}= Int(a_\gamma^{-1}\gamma(n_w)a_\gamma),$

{\bf (iii)} $\wt\mu(\gamma)(n_w)\equiv a_\gamma^{-1}\gamma(n_w)a_\gamma\quad\in N_{\bG}(\bT)/\bT,$

\noindent where $a_\gamma$ as in {\bf (3a), (3b)}.
\endproclaim
\demo{Proof} {\bf (iii)} is a consequence of (i), (ii) because $Int: N_\bG(\bT)/\bT \rightarrow Aut\;X^*(\bT)$ is injective.

\noindent {\bf (i)} For $\chi\in X^*(\bT)$ we use that $\mu(\gamma)(\chi)=\chi\circ\wt\mu(\gamma)^{-1}.$ Then:
$$ \mu(\gamma)\cdot w\cdot\mu(\gamma)^{-1}\circ\chi =(w\cdot(\chi\circ\wt\mu(\gamma)))\circ\wt\mu(\gamma^{-1})= n_w(\chi\circ\wt\mu(\gamma))n_w^{-1}\circ\wt\mu(\gamma^{-1}).$$
On the other hand
$$n_w(\chi\circ\wt\mu(\gamma))n_w^{-1}(\wt\mu(\gamma^{-1})\circ t)=(\chi\circ\wt\mu(\gamma))(n_w^{-1}(\wt\mu(\gamma^{-1})\circ t)n_w)=$$
$$=\chi(\wt\mu(\gamma)(n_w^{-1})\cdot t\cdot\wt\mu(\gamma)(n_w))=
((\wt\mu(\gamma)(n_w))\chi(\wt\mu(\gamma)(n_w))^{-1})(t),$$
which yields the proof of (i).

\noindent {\bf (ii)} The operator $\mu(\gamma)^{-1}=\mu(\gamma^{-1})\in Aut(X^*(\bT))$ is inverse to $\chi\mapsto\gamma\chi\gamma^{-1}\circ Int(a_\gamma).$ By a direct computation we see:
$$  \mu(\gamma)^{-1}=\{\chi\mapsto \gamma^{-1}\circ\chi(a_\gamma^{-1}\gamma(*)a_\gamma)\}.$$
Therefore:
$$ \mu(\gamma)\cdot w\cdot\mu(\gamma)^{-1}\circ\chi=\mu(\gamma)\cdot w(\gamma^{-1}\circ\chi(a_\gamma^{-1}\gamma(*)a_\gamma))=$$
$$=\mu(\gamma)\circ (\gamma^{-1}\circ\chi(a_\gamma^{-1}\gamma(n_w^{-1}*n_w)a_\gamma))=
\mu(\gamma)\circ(\gamma^{-1}\circ\chi(a_\gamma^{-1}\gamma(n_w)^{-1}\gamma(*)\gamma(n_w)a_\gamma))=$$
$$ =\chi(a_\gamma^{-1}\gamma(n_w)^{-1}\cdot\gamma(\gamma^{-1}(a_\gamma*a_\gamma^{-1}))\cdot\gamma(n_w)a_\gamma)=\chi(a_\gamma^{-1}\gamma(n_w)^{-1}a_\gamma\cdot *\cdot a_\gamma^{-1}\gamma(n_w)a_\gamma),$$
which proves (ii).\qed\enddemo
\proclaim{2.2.3 Corollary} We may turn $W(\bT,\bG)\cong W(\wh T,\wh G)$ into an isomorphism
$$  N_\bG(\bT)/\bT \cong N_{\wh G}(\wh T)/\wh T,\qquad w\mapsto \wh w,$$
which is $\wt\mu(\Gamma)$-$\Gamma$-equivariant, i.e.
$\wh{\wt\mu(\gamma)(w)} = \gamma(\wh w)$,

\noindent with the notation $\gamma(\wh w):= \wt\mu_{\wh G}(\gamma)(\wh w).$
\endproclaim
\demo{Proof}We note that 2.2.2(i) also makes sense with respect to $W(\wh T,\wh G)\cong N_{\wh G}(\wh T)/\wh T,$ $\wh w\mapsto n_{\wh w}$ and yields the equality:
$$  \mu_{\wh G}(\gamma)\wh w \mu_{\wh G}(\gamma)^{-1} = Int(\wt\mu_{\wh G}(\gamma)(n_{\wh w})).$$
Therefore the Galois-equivariance of {\bf (4)} turns into our assertion where we have replaced now $n_w,$ $n_{\wh w}$ by the notations $w$ and $\wh w$ resp. \qed\enddemo

\bigskip
{\bf 2.3 Parabolic subgroups and Levi subgroups of $G$ and $\wh G$ resp.}

\noindent Since $(\bG,\;\bB,\; N_\bG(\bT),\;\Delta^*(\bT,\bB))$ is a BN-pair, where we identify the simple roots with the corresponding reflections in $W=W(\bT,\bG),$ we obtain an inclusion preserving bijection
$$   wW_I \leftrightarrow \bP(wW_I):=w\bB W_I\bB w^{-1}$$
between cosets in $W$, where $W_I$ corresponds to the subset $I\subseteq \Delta^*(\bT,\bG),$ and semistandard parabolic subgroups $\bP\supset\bT.$ Therefore the isomorphism {\bf (4)} and its variant {\bf 2.2.3} induce now a bijection
$$ \bT\subset\bP\leftrightarrow wW_I\longleftrightarrow \wh w\wh W_{\wh I}\leftrightarrow
\wh P\supset \wh T\tag{5}$$ between semi-standard parabolic groups
$\bP\supset\bT$ and semi-standard groups $\wh P\supset\wh T,$ such that $\bP(wW_I)$ corresponds to $\wh P(\wh w\wh W_{\wh I})=\wh w\wh B\wh W_{\wh I}\wh B\wh w^{-1}$ where $\wh I\leftrightarrow I$ are corresponding subsets of $\Delta^\ast(\wh T,\wh B)$ and $\Delta^\ast(\bT,\bB)$ resp.. (Note here that the notation $\bP(wW_I)=w\bB W_I\bB w^{-1}$ refers to $W(\bT,\bG)=N_\bG(\bT)/\bT,$ therefore $w$ means now $n_w$; and similar for $\wh G.$)
Moreover the semi-standard Levi group $\bM_\bP$ of $\bP(wW_I)$ is characterized by the fact that
$$ W(\bT,\bM_\bP) = wW_Iw^{-1},$$
and we will write this as $\bM_\bP=\bM(wW_Iw^{-1}).$ So our bijection of semi-standard parabolics induces a bijection
$$  \bM= \bM(wW_Iw^{-1})\leftrightarrow \wh M=  \wh M(\wh w\wh W_{\wh I}\wh w^{-1})\tag{5a}$$
of semi-standard Levi groups. We note here that $\wh M$ is actually dual to $\bM$ in the sense that
$$  \psi_0(\bM,\bM\cap w\bB w^{-1},\bT)=\psi_0(\wh M,\wh M\cap \wh w\wh B\wh w^{-1},\wh T),$$
where $\Delta^*(\bT,\bM\cap w\bB w^{-1})= wI \subset\Delta^*(\bT,w\bB w^{-1})$ is dual to $\Delta^*(\wh T,\wh M\cap \wh w\wh B \wh w^{-1})= \wh w\wh I \subset\Delta^*(\wh T,\wh w\wh B \wh w^{-1}).$

\noindent {\bf Remark:} A semi-standard Levi group is recovered from its Weyl group via
$\bM = C_\bG((\bT^{W(\bT,\bM)})^0).$ Therefore the bijection {\bf (5a)} between semi-standard Levi groups is characterized by the fact that under the distinguished isomorphism $W(\bT,\bG)\cong W(\wh T,\wh G)$ (which is based on the identification $\Delta^*(\bT,\bB)\leftrightarrow\Delta(\wh T,\wh B)$) we have
$$   W(\bT,\bM) \cong W(\wh T,\wh M).$$
In particular $\bM$ is standard if and only if $\wh M$ is.

\noindent Of course these bijections are equivariant with respect to the Weyl-group actions
$$  w'\bP(wW_I){w'}^{-1} =\bP(w'w W_I),\qquad w'\bM(wW_I w^{-1}){w'}^{-1} =\bM(w'wW_I w^{-1} {w'}^{-1})$$
on the $\bG$-side and on the $\wh G$-side resp. But from {\bf 2.2.3} we see
\proclaim{2.3.0} that the bijections (5), (5a) are also equivariant with respect to the $\wmu(\Gamma)$- and $\Gamma$-action resp., i.e.
$$  \;^{\wmu(\gamma)}\bP(wW_I) = \bP(\;^{\wmu(\gamma)}(wW_I))\leftrightarrow \;^\gamma{\wh P}(\wh w\wh W_{\wh I})=\wh P(\;^\gamma(\wh w\wh W_{\wh I})),$$
and similar for the Levi groups.\endproclaim

\noindent The standard groups $\bP\supseteq\bB$ are in bijection with the subsets $I\subseteq \Delta^*(\bT,\bB):$
$$  I\mapsto \bP_I=\bP(W_I),\qquad  \bP\mapsto I(\bP)=\Delta^*(\bT,\bB) -\Delta^*(\bT,\bP),$$
where the last difference means to remove those simple roots which come from the unipotent radical of $\bP.$ Similar on the $\wh G$-side. Finally we still remark that the bijections are inclusion preserving:

{\bf (i)} $w_1W_{I_1}\subseteq w_2W_{I_2}$ if and only if $\bP(w_1W_{I_1})\subseteq\bP(w_2W_{I_2}),$

{\bf (ii)} $w_1W_{I_1}w_1^{-1}\subseteq w_2W_{I_2}w_2^{-1}$ if and only if $\bM(w_1W_{I_1}w_1^{-1})\subseteq \bM(w_2W_{I_2}w_2^{-1}).$

In (i) from right to left one uses that a parabolic group containing a standard group is standard by itself, and in (ii) one uses that $\bM_1\subseteq\bM_2$ if and only if $W(\bT,\bM_1)\subseteq W(\bT,\bM_2).$

On $(\bG,\bB,\bT)$ we have to distinguish the usual $\Gamma$-action and the $\wmu(\Gamma)$-action. In order to compare them we will need:
\proclaim{2.3.1 Lemma} Consider $\gamma\in \Gamma$ and $a_\gamma\in \bG$ such that $a_\gamma(\bB,\bT)a_\gamma^{-1}= (\gamma\circ\bB,\gamma\circ\bT).$ Then for the semi-standard groups $\bP=\bP(wW_I),$ $\bM_\bP= \bM(wW_Iw^{-1})$ we obtain
$$  a_\gamma(\bP,\bM_\bP)a_\gamma^{-1}=  (\gamma\circ\wmu(\gamma)^{-1}\circ\bP,\;\gamma\circ\wmu(\gamma)^{-1}\circ\bM_\bP).$$
\endproclaim
\demo{Proof} Here we consider $W(\bT,\bG)=N_\bG(\bT)/\bT,$ that means we identify $w$ and $n_w.$ Then according to {\bf 2.2.2} the Galois action rewrites as
$$  \wmu(\gamma)\circ w= a_\gamma^{-1}\gamma(w) a_\gamma,\qquad a_\gamma w a_\gamma^{-1}=\gamma(\wmu(\gamma)^{-1}\circ w),$$  hence
$$  a_\gamma\bP(wW_I)a_\gamma^{-1}=a_\gamma w\bB W_I\bB w^{-1}a_\gamma^{-1}=$$
$$=a_\gamma wa_\gamma^{-1}\cdot a_\gamma\bB a_\gamma^{-1}\cdot a_\gamma W_I a_\gamma^{-1}\cdot a_\gamma\bB a_\gamma^{-1}\cdot(a_\gamma wa_\gamma^{-1})^{-1}=$$
$$=\gamma(\wmu(\gamma)^{-1}\circ w)\cdot \gamma(\bB)\cdot\gamma(\wmu(\gamma)^{-1}\circ W_I)\cdot\gamma(\bB)\cdot\gamma(\wmu(\gamma)^{-1}\circ w)^{-1}=\gamma\circ\wt\mu(\gamma)^{-1}\circ\bP(wW_I),$$
where we have used that $\gamma(\bB)= \gamma\circ\wt\mu(\gamma)^{-1}\circ\bB.$

Finally $a_\gamma\bM_\bP a_\gamma^{-1}$ is the uniquely determined semistandard Levi-complement of $a_\gamma\bP a_\gamma^{-1}$ with respect to the maximal torus $a_\gamma\bT a_\gamma^{-1}=\gamma\circ\wt\mu(\gamma)^{-1}\circ\bT$, and therefore the relation for $a_\gamma\bM_\bP a_\gamma^{-1}$ follows.
\qed\enddemo

\noindent {\bf 2.3.2} For the rest of the paper we choose the minimal pair of $F$-groups $(\bP_0,\bM_0)$ from {\bf 1.1} and the pair $(\bB,\bT)$ from {2.1} in a compatible way, that means the groups $(\bP_0,\bM_0)$ should be standard with respect to $(\bB,\bT)$. More precisely we fix

\noindent $(\bP_0,\bM_0,\bA_0)$ as in {\bf 1.1}. Then we choose a maximal torus $\bT\subseteq \bM_0$ and a corresponding Borel-subgroup $\bB$ in $\bP_0,$ hence
$$   \bP_0\supseteq\bB\supset \bT\supseteq\bA_0.\tag{*}$$
Note that $\bA_0$ is the unique maximal $F$-split subtorus in $\bM_0$ (because it is central), and therefore it is also the unique maximal $F$-split subtorus in $\bT.$
\proclaim{Lemma} If $(\bB',\bT')$ is another pair with property {\bf (*)}, then $(\bB',\bT')=x(\bB,\bT)x^{-1}$ for some $x\in \bM_0,$ where the coset $x\bT$ is uniquely determined. In particular we may take $(\bB',\bT'):=(\gamma(\bB),\gamma(\bT)),$ for any $\gamma\in \Gamma$
and therefore in {\bf (3a),(3b)} we will always have $a_\gamma\in\bM_0.$
\endproclaim
\demo{Proof} We are given two pairs consisting of a
Borel subgroup and a maximal torus which are both in $\bP_0.$ Therefore according
to {\bf [Sp2],6.4.1 and 6.4.12} we may find $x\in\bP_0.$ Furthermore $\bM_0\supset \bT'$ implies now $x^{-1}\bM_0 x\supset\bT$ hence $x^{-1}\bM_0x=\bM_0$ because in $\bP_0$ there is only one semi-standard Levi-subgroup. Finally writing $x=m_0u_0\in\bP_0,$ where $u_0$ is from the unipotent radical, we see from {\bf [DM], 1.18} that $x^{-1}\bM_0 x=\bM_0$ implies $u_0=1,$ hence $x=m_0\in\bM_0.$ If $(\bB',\bT')=y(\bB,\bT)y^{-1}$ is any solution then always $y\in x\bT\subseteq\bM_0.$
\qed\enddemo

\proclaim{2.3.3 Corollary} Let $\bM\subseteq\bG|F$ be a semi-standard $F$-Levi-subgroup of $\bG$
and consider the injection  $X^*(\bM)\rightarrow X^*(\bT)$ which is given via restriction.
Then the $\mu(\Gamma)$-action on $X^*(\bT)$
and the usual Galois-action on $X^*(\bM)$ agree.\endproclaim
\demo{Proof} According to {\bf (3a)} we have
$$ (\gamma\cdot_\mu \chi)(t) = (\gamma\chi\gamma^{-1})(a_\gamma ta_\gamma^{-1})=\gamma\circ\chi(\gamma^{-1}(a_\gamma)\cdot\gamma^{-1}(t)\cdot\gamma^{-1}(a_\gamma^{-1})),$$
where $a_\gamma\in\bM_0\subseteq\bM.$ Therefore $\chi\in X^*(\bM)$ will imply
$  \gamma\cdot_\mu \chi = \gamma\chi\gamma^{-1}.$\qed\enddemo

Moreover we obtain now the following characterization of
those semi-standard groups $\bP,\;\bM\supset\bT$ which are defined over $F:$

\proclaim{2.3.4 Proposition} Let $\bP,\;\bM\supseteq\bT$ be semi-standard parabolic and Levi-subgroups resp. Then $\bP,$ $\bM$ resp. is defined over $F$ if and only if it contains $\bM_0$ and is stable under the $\wmu(\Gamma)$-action.\endproclaim
\demo{Proof} We argue with $\bP.$ According to {\bf 2.3.1} we have
$$  a_\gamma\bP a_\gamma^{-1}= \gamma\circ\wmu(\gamma)^{-1}\circ\bP.\tag{*}$$
for any $\gamma\in \Gamma.$
If $\bP$ is a semi-standard $F$-group then it contains $\bM_0$ and therefore by the Lemma of 2.3.2  the displayed equality rewrites as $\bP= \gamma\circ\wmu(\gamma)^{-1}\circ\bP$ for all $\gamma\in\Gamma.$ Therefore $\bP$ must be $\wmu(\Gamma)$-stable because it is $\Gamma$-stable.

\noindent Conversely assume that $\bP\supset\bM_0$ and that it is stable under the $\wt\mu(\Gamma)$-action. Then it is obvious from {\bf (*)} that it is also stable under the $\Gamma$-action, hence it is defined over $F.$

\noindent The same argument works for the semi-standard Levi-subgroups because again we have:
$ a_\gamma\bM a_\gamma^{-1}= \gamma\circ\wmu(\gamma)^{-1}\circ\bM.$
\qed\enddemo
\proclaim{2.3.5 Corollary} Let $\wh P,\;\wh M\supset\wh T$ be semi-standard groups in $\wh G.$ Then under (5), (5a) resp. they correspond to semistandard groups $\bP,\;\bM\supset\bT$ which are defined over $F$ if and only if:

\noindent $\wh P,\;\wh M$ are $\Gamma$-stable and contain $\wh M_0.$\endproclaim
\demo{Proof} This follows directly from 2.3.4 and 2.3.0.\enddemo
\bigskip

{\bf 2.4  $F$-parabolic groups and their $\wh G$-counterparts in terms of cosets in $W,$ $\wh W$ resp.}

\proclaim{2.4.1 Lemma} {\bf (i)} A semi-standard parabolic $\bP=\bP(wW_I)$ is $\wmu(\Gamma)$-stable if and only if $wW_I$ is $\wmu(\Gamma)$-stable, or equivalently:

\noindent $I\subset \Delta^*(\bT,\bB)$ is a $\mu(\Gamma)$-stable subset and $wW_I=w'W_I$ for some $w'\in W^{\wmu(\Gamma)}$.

\noindent {\bf (ii)} A semi-standard Levi group $\bM$ is $\wt\mu(\Gamma)$-stable if and only if it comes as $\bM=\bM_\bP=\bM(wW_I w^{-1})$ for a semi-standard parabolic $\bP=\bP(wW_I)$ which is $\wt\mu(\Gamma)$-stable.
\endproclaim
\demo{Proof} {\bf (i)} The first statement is obvious and concerning the second statement we only need to show that the condition is necessary. What is obvious (from 2.2.2(i)) is that $I$ has to be $\mu(\Gamma)$-stable if the coset $wW_I$ is $\wmu(\Gamma)$-stable. So we are left to show that we can find a $\wmu(\Gamma)$-invariant representative in $wW_I$.

For this we begin with the case where $(\bG,\bB,\bT)$ is quasisplit, hence $\bM_0=\bT.$ Then from {\bf 2.3.3} we see that on $X^*(\bT)$  the usual $\Gamma$-action and the $\mu(\Gamma)$-action agree. This takes over from $X^*(\bT)$ to $W$ because of 2.2.2(iii) together with $a_\gamma\in \bM_0=\bT.$
But then it follows that $\bP=\bP(wW_I)$ is $\Gamma$-stable, {\bf hence an $F$-parabolic group,} if and only if it is $\wmu(\Gamma)$-stable. And an $F$-parabolic group will contain an $F$-Borelgroup, hence
$$  \bP(wW_I) \supseteq \bP(w')=w'\bB{w'}^{-1},$$
and therefore $w'\in wW_I$ has to be $\wmu(\Gamma)$-stable.

Now we go to the general case. Then $(\bG,\bB,\bT)$ is an inner form of a uniquely determined quasisplit group $(\bG',\bB',\bT'),$ in particular we have then an $\overline F$-isomorphism $f:(\bG,\bB,\bT)\rightarrow (\bG',\bB',\bT')$  and according to {\bf 2.2.1} the induced isomorphism
 $$  w\in W(\bT,\bG)\mapsto (f^*)^{-1}\cdot w\cdot f^*\in W(\bT',\bG')$$
is $\mu(\Gamma)$-equivariant. Therefore we can transport our question to the quasisplit side where it has been answered already.

\noindent {\bf (ii)} Here we argue in the same way as in {\bf (i)}. If $\bG$ is quasisplit, then $\bP$, $\bM$ are $\wt\mu(\Gamma)$-stable if and only if they are $\Gamma$-stable, and (ii) is true if we replace $\wt\mu(\Gamma)$-stable by the property $F$-parabolic. Then in general, $\bG$ is the inner form of a quasisplit group $\bG'$ and $\bM=\bM(wW_Iw^{-1})$ is $\wmu(\Gamma)$-stable if and only if the Weyl group $wW_Iw^{-1}$ has this property. Therefore again we may use {\bf 2.2.1} to shift our assertion to the quasisplit case.
\qed\enddemo

The next step is to find an appropriate expression of the relative Weyl group
$$\;_FW(\bA_0,\bG):=N_\bG(\bA_0)/C_\bG(\bA_0)= N_G(A_0)/C_G(A_0),\tag{7}$$
where $G=\bG(F)$ and $A_0=\bA_0(F),$ in terms of the absolute Weyl group $W=W(\bT,\bG).$
(See {\bf [Bo2],21.2} for the
fact that elements of the relative Weyl group can be represented by rational points.)
First of all our setting $\bP_0\supseteq\bB\supset\bT\supseteq\bA_0$ implies that $\bP_0,$ $\bM_0$ are given as standard groups
$$  \bP_0=\bP(W_{I_0}),\qquad  \bM_0=\bM(W_{I_0}),$$
where $I_0=\Delta^*(\bT,\bB\cap\bM_0)\subset \Delta^*(\bT,\bB)$ is the subset of simple roots in $\bM_0,$ which has to be $\mu(\Gamma)$-stable as we see from the previous Proposition.

\proclaim{2.4.2 Proposition}  Let $W_0:=W_{I_0}=W(\bT,\bM_0)\subseteq W=W(\bT,\bG)$ be the absolute Weyl groups and consider the normalizer $N_W(W_0) = N_W(\bM_0).$ Then the relative Weyl group $\;_FW$ can be written as
$$  \;_FW \cong (N_W(W_0)/W_0)^{\mu(\Gamma)}\cong N_W(W_0)^{\mu(\Gamma)}/W_0^{\mu(\Gamma)},$$
where the first isomorphism is the map $N_\bG(\bA_0)/C_\bG(\bA_0)@<\sim <<  N_W(\bA_0)/C_W(\bA_0),$ and the second isomorphism comes from {\bf 2.4.1(i)}.
\endproclaim

\demo{Sketch of Proof} We have a commutative diagram

$$ \CD
N_W(\bA_0)/C_W(\bA_0)@>>> N_W(\bM_0)/C_W(\bA_0)\\
@VV\cong V               @VV\cong V\\
\;_FW=N_\bG(\bA_0)/C_\bG(\bA_0)@>>> N_\bG(\bM_0)/C_\bG(\bA_0)
\endCD $$
where the horizontal arrows are {\bf inclusions}. And in the upper row we have:
$$  C_W(\bA_0)=W(\bT,\bM_0)= W_0,\qquad N_W(\bM_0)=N_W(W_0).$$
Therefore the right vertical of the diagram can be written as
$$    N_\bG(\bM_0)/\bM_0 \cong N_W(W_0)/W_0 ,$$
which is $\Gamma$- $\mu(\Gamma)$-equivariant, because from the proof of {\bf 2.3.1} and the Lemma of {\bf 2.3.2} we see that $\mu(\gamma)\circ w= a_\gamma^{-1}\cdot\gamma(w)\cdot a_\gamma\equiv \gamma(w)\;\text{mod}\;\bM_0.$ So we see that $$(N_\bG(\bM_0)/\bM_0)^\Gamma\cong (N_W(W_0)/W_0)^{\mu(\Gamma)}\cong N_W(W_0)^{\mu(\Gamma)}/W_0^{\mu(\Gamma)}$$
where the last isomorphism follows from 2.4.1(i). On the other hand we have
$$  N_G(A_0)=N_\bG(\bA_0)^\Gamma= N_\bG(\bM_0)^\Gamma= N_G(M_0)$$
and therefore from (7) we see that under the lower horizontal of our diagram  $\;_FW$ is mapped onto the subgroup of $N_\bG(\bM_0)/\bM_0$ which is generated by $N_\bG(\bM_0)^\Gamma.$  To see that the inclusion $N_\bG(\bM_0)^\Gamma\cdot\bM_0/\bM_0\subseteq (N_\bG(\bM_0)/\bM_0)^\Gamma$ is actually an equality we consider the set $\PPP(\bM_0)$ of all parabolic groups with Levi subgroup $\bM_0$ and the injection
$$ N_\bG(\bM_0)/\bM_0\hookrightarrow \PPP(\bM_0),\qquad g\mapsto g\bP_0 g^{-1}.$$
Then the assertion follows from the fact that $\;_FW$ acts simply transitive on the subset of $F$-parabolic groups from $\PPP(\bM_0).$
\qed\enddemo

\proclaim{2.4.3 Corollary}  A semi-standard parabolic group $\bP=\bP(wW_I) \supset\bT$ is $F$-parabolic if and only if the following conditions hold:

{\bf (i)} the subset $I\subseteq \Delta^*(\bT,\bB)$ is stable under the action of $\mu(\Gamma),$ and $I\supseteq I_0,$  equivalently $W_I\supseteq W_0$ is $\wmu(\Gamma)$-stable,

{\bf (ii)} $   wW_I \cap  N_W(W_0)^{\wmu(\Gamma)} \ne\emptyset ,$
this means we have a representative $w'\in wW_I$ which normalizes $W_0,$ and which is fixed under $\wmu(\Gamma),$ {\bf hence:} $wW_I \cap  N_W(W_0)^{\wmu(\Gamma)}=w'\cdot N_{W_I}(W_0)^{\wmu(\Gamma)}.$
\endproclaim
\demo{Proof} This is now an easy consequence of {\bf 2.3.4, 2.4.1, 2.4.2} and the well known
\proclaim{Lemma} Every semi-standard $F$-parabolic group $\bP\supset\bM_0$ occurs as an $\,_FW$-conjugate of a uniquely determined standard $F$-parabolic $\bP'\supseteq \bP_0.$\qed\endproclaim
\enddemo
The cosets $wW_I\subseteq W$ which fulfill the conditions of 2.4.3 we may call relevant cosets. It is now immediate to define the notion of a relevant coset $\wh w\wh W_{\wh I}\subseteq \wh W$ and to express the groups $\wh P,\;\wh M$ of {\bf 2.3.5} in terms of relevant cosets. In the case of Levi groups one uses {\bf 2.4.1(ii)}.
\bigskip

{\bf 2.5 Parabolic subgroups and Levi subgroups of $\;^LG$}

We recall the definitions of {\bf [Bo1], 3.3, 3.4} and add some remarks which come
from the previous subsections.

\noindent A parabolic subgroup $P(\;^LG)$ is by definition a closed
subgroup $P(\;^LG)\subset \;^LG$ of the form $P(\wh
G)\hookrightarrow P(\;^LG)\twoheadrightarrow W_F,$ where $P(\wh
G)\subset \wh G$ is a parabolic subgroup which has a maximal
normalizer in $\;^LG.$ (Note that for any parabolic $P(\wh G)$ with
normalizer $N$ in $\;^LG$ we must have $N/P(\wh G)\hookrightarrow
\;^LG/\wh G \cong W_F.$)

\noindent It is called {\bf standard (semi-standard)} if
$P(\;^LG)\supseteq \;^LB=\wh B\rtimes W_F$, or $P(\;^LG)\supset
\;^LT=\wh T\rtimes W_F$ resp., which implies $P(\;^LG)=P(\wh G)\rtimes
W_F$ such that $P(\wh G)$ is standard (semi-standard) in $(\wh G,\wh
B,\wh T)$ and is stable under the $\Gamma$-action.

\noindent Note here that $P(\;^LG)\supseteq \;^LB=\wh B\rtimes
W_F$ implies that the connected component of $P(\;^LG)$ which is a
parabolic subgroup of $\wh G$ must contain the connected component
of $\;^LB$ which is $\wh B$, and similar for $\;^LT$. Therefore:
\proclaim{2.5.1 Proposition} {\bf (i)} The bijection {\bf (5)}:
$\bT\subset\bP\leftrightarrow \wh P\supset\wh T$ induces a
bijection
$$ \bP\leftrightarrow \;^LP=\wh P\rtimes W_F$$
between subgroups $\bP\supset \bT$ which are $\wmu(\Gamma)$-stable
and semi-standard parabolic groups
$P(\;^LG)\supset\;^LT.$

\noindent {\bf (ii)} The standard parabolic groups in $\;^LG$ are
precisely the groups
$$P(\;^LG)=P(\wh G)\rtimes W_F,$$ where $P(\wh G)\supseteq \wh B$ is
standard and such that $I(P(\wh G))\subseteq \Delta^*(\wh T,\wh
B)=\Delta_*(T,B)$ is stable under the $\Gamma$-action.\endproclaim

\noindent This is
because $\wh B$ itself is stable under $\Gamma$-action, hence a
standard parabolic $P(\wh G)$ is stable under $\Gamma$-action if and
only if $I(P(\wh G))$ is.

And from {\bf [Bo1] 3.3} we know that every parabolic subgroup of $\;^LG$ is $\wh G$-conjugate to precisely one standard parabolic subgroup.
\bigskip

We recall that the
notion of unipotent radical and of Levi complement is defined for
any algebraic group which need not be connected. (See e.g. {\bf
[DM]0.16(ii) and 1.16 resp.}) \noindent {\bf A Levi subgroup}
$L(\;^LG)\subset \;^LG$ is by definition a Levi complement of the
unipotent radical of some parabolic subgroup $P(\;^LG)$ as
considered before.
\proclaim{2.5.2 Remark} If the parabolic group $P(\;^LG)$ is
semi-standard (i.e. if $P(\;^LG)\supset\;^LT$) then it contains
precisely one Levi complement $L(\;^LG)\supset \;^LT,$ namely
$$  L(\;^LG) = L(\;^LG)^0\rtimes W_F,$$
where $L(\;^LG)^0=L(\wh G)$ is the uniquely determined Levi
subgroup of $P(\;^LG)^0= P(\wh G)$ which contains $\wh
T.$\qed\endproclaim

The Levigroup $L(\;^LG)$ is called {\bf semi-standard / standard} if it is the distinguished Levi complement in a semi-standard / standard
parabolic  group $P(\;^LG)$ which means we must have
$$   P(\;^LG)\supseteq L(\;^LG)\supseteq \wh T\rtimes W_F.$$
\proclaim{2.5.3 Remark} A Levi subgroup $L(\;^LG)$ is always $\wh G$-conjugate to at least one standard Levi subgroup $L'(\;^LG)= L'(\wh G)\rtimes W_F.$\qed\endproclaim

Then from {\bf Proposition 2.3.4} we see:

\proclaim{2.5.4 Proposition} {\bf (i)} The bijections {\bf (5):} $\bT\subset\bP\leftrightarrow \wh P\supset \wh T$ and {\bf (5a):}
$\bT\subset\bM\leftrightarrow \wh M\supset \wh T$ between semi-standard parabolic and Levi-subgroups in $\bG$ and $\wh G$ resp. induce bijections
$$  \bP\leftrightarrow \;^LP=\wh P\rtimes W_F,\qquad \bM\leftrightarrow \;^LM=\wh M\rtimes W_F,\tag{F}$$
between the semi-standard $F$-groups on one hand  and those semi-standard subgroups of $\;^LG$ which contain $\;^LM_0$ on the other hand.

{\bf (ii)} In the same way as the bijections (5), (5a) are equivariant with respect to the action of $W(\bT,\bG)\cong W(\wh T,\wh G),$ the bijections {\bf (F)} are now equivariant with respect to the action of $\;_FW=N_W(W_0)^{\wmu(\Gamma)}/W_0^{\wmu(\Gamma)}\cong \;_F\wh W:=N_{\wh W}(\wh W_0)^\Gamma/\wh W_0^\Gamma.$

\endproclaim

For the proof one has only to translate what it means on the $\;^LG$-side that a semi-standard group $\bP$ or $\bM$ is $\wmu(\Gamma)$-stable and contains $\bM_0.$

\noindent The notation $\;^LM$ is justified because the identification $\psi_0(\bM,\bM\cap w\bB w^{-1},\bT)=\psi_0(\wh M,\wh M\cap \wh w\wh B\wh w^{-1},\wh T)$ is $\wt\mu(\Gamma)$-$\Gamma$ equivariant if $\bM=\bM(wW_Iw^{-1})$ is an $F$-group.

\definition{2.5.5 Definition and Remark} The semi-standard parabolic and Levi subgroups of $\;^LG$ which come from semi-standard $F$-subgroups in $\bG$ and which are characterized by {\bf 2.5.4}, we call the {\bf relevant} semistandard subgroups of $\;^LG.$ Then a parabolic or Levi subgroup of $\,^LG$ will be relevant in the sense of {\bf [Bo1], 3.3, 3.4} if and only if it is $\wh G$-conjugate to a relevant semi-standard group.
\enddefinition

The connection between our relevant semi-standard groups and the relevant groups in the sense of loc.cit. is an immediate consequence of
\proclaim{2.5.6 Lemma} {\bf (i)} For a parabolic group $P(\,^LG)$ the following is equivalent:

{\bf a)} There exists $x\in \wh G$ such that $xP(\,^LG)x^{-1}\supseteq \,^LP_0,$

{\bf b)} There exists $x\in \wh G$ such that $xP(\,^LG)x^{-1}\supseteq \,^LM_0.$

\noindent {\bf (ii)} $L(\,^LG)$ occurs as Levi-subgroup of a parabolic $P(\,^LG)$ such as in {\bf (i)}, if and only if there exists $x\in \wh G$ such that $xL(\,^LG)x^{-1}\supseteq \,^LM_0.$
\endproclaim
\demo{Proof} In (i) we only need to prove that b) implies a). In view of {\bf 2.5.4} the assumption b) means: $xP(\,^LG)x^{-1} = \,^LP$ where $\bP$ is a semi-standard $F$-parabolic group. But then $\bP$ is $\,_FW$-conjugate to a standard $F$-parabolic and therefore $\,^LP$ is $\,_F\wh W$-conjugate to a group which contains $\,^LP_0.$

\noindent In (ii) it is clear that $xL(\,^LG)x^{-1}\supseteq \,^LM_0$ implies $xP(\,^LG)x^{-1}\supset \,^LM_0.$ Conversely assume that $xP(\,^LG)x^{-1}\supset \,^LM_0$ and that $L(\,^LG)$ is a Levi subgroup of $P(\,^LG).$ Then, using {\bf 2.5.4}, we see that $xP(\,^LG)x^{-1}$ must have a Levi-subgroup $\,^LM\supseteq \,^LM_0$ and therefore $L':= x^{-1}\,^LM x$ is a Levi subgroup of $P(\,^LG)$ which is good. But any other Levi subgroup of $P(\,^LG)$ must be $\wh G$-conjugate to $L'.$ Here we may use the corresponding fact for the subgroups in $\wh G$ (see for instance {\bf [DM], 1.18}) and the fact that $L(\,^LG)= N_{P(\,^LG)}(L(\wh G))$ comes as the normalizer of a usual Levi subgroup, if $L(\,^LG)$ is a Levi subgroup of $P(\,^LG).$ As a consequence we see that together with $L'$ also any other Levi subgroup of $P(\,^LG)$ will be good.\qed\enddemo

\bigskip
{\bf 3. L-parameters and their Langlands classification}

We repeat the basic definitions:
\definition{3.1 Definition} {\bf (i)} An {\bf $L$-homomorphism} is a homomorphic map
$$  \phi=(\phi_1,\phi_2): SL_2(\C)\times W_F \rightarrow \;^LG=\wh G\rtimes W_F$$
from the direct product (therefore two separate maps) to the semidirect product, such that $Im(\phi_1),$ $Im(\phi_2)$ commute with each other and moreover:

$\phi_1:SL_2(\C)\rightarrow \wh G$ is a rational morphism of complex algebraic groups (hence the image in $\;^LG$ must be connected),

$\phi_2:W_F\rightarrow \wh G\rtimes W_F$ is a continuous homomorphism of the form
$$  \phi_2(w)= (\vf_2(w),w)$$
which means that $\vf_2:W_F\rightarrow \wh G$ is a 1-cocycle, and $\phi_2(w)\in\;^LG$ should be always a
{\bf semisimple element} i.e. under any rational representation $r:\;^LG\rightarrow GL_N(\C),$ it is mapped to a semisimple element. (One can restrict here to the case where $w\in W_F$ is a Frobenius-lift and one speaks of Frobenius-semisimple maps $\phi$ and $\phi_2$ resp.)

\noindent Finally, if $\bG$ is not quasisplit and $Im(\phi)\subset L(\;^LG)$ then $L(\;^LG)$ should be a {\bf relevant} Levi subgroup in $\;^LG.$ (see 2.5.5 above).

{\bf 3.1 (ii)} On the set of $L$-homomorphisms we have an obvious {\bf action of $\wh G$ by conjugation.} If $x\in \wh G$ and $\phi$ are given then we define the {\bf equivalent} $L$-homomorphism $\;^x\phi$ as:
$$ \;^x\phi(A,w):= x\phi_1(A)\phi_2(w)x^{-1}=x\phi_1(A)x^{-1}x\phi_2(w)x^{-1}=\;^x\phi_1(A)\;^x\phi_2(w),$$
for all $A\in SL_2(\C),$ $w\in W_F.$ Obviously $\;^x\phi =\phi$ if and only if $x\in C_{\wh G}(Im(\phi))$ is in the $\wh G$-centralizer of the image of $\phi.$

{\bf 3.1 (iii)} One may also use the restriction of $\phi$ to the Weil-Deligne group which is obtained via the embedding
$$  \C^+\rtimes W_F\hookrightarrow SL_2(\C)\times W_F,$$
$$ z\in\C \mapsto \bmatrix 1&z\\ 0&1\endbmatrix \times 1_W,\qquad w\in W_F\mapsto w_*:= \bmatrix \pa w\pa^{1/2}& 0\\0&\pa w\pa^{-1/2}\endbmatrix \times w,$$
where $\pa w\pa= q^{d(w)}$ if $\overline w\in W_F/I_F$ is the $d(w)$-th power of the Frobenius automorphism and $q$ is the order of the residue field of $F.$ (see {\bf [L1],p.209} and {\bf [T],(1.4.1),(1.4.6)} resp.).
We will denote  $ W_{F*}:=\{w_*\;|\; w\in W_F\}$
the second copy of $W_F$ in $SL_2(\C)\times W_F$ which comes from the Weil-Deligne group. If now $\phi=(\phi_1,\phi_2)$ is an L-homomorphism, then we have
$$   \phi(w_*)=\phi_1\left(\bmatrix \pa w\pa^{1/2}&0\\ 0&\pa w\pa^{-1/2}\endbmatrix\right)\cdot \phi_2(w)\quad\in \;^LG$$
the product of two commuting semisimple elements, hence $\phi(w_*)$ is again semisimple.\enddefinition

\definition{Remark} In {\bf 3.1 (i)} it is possible to replace $\,^LG$ by $\wh G\rtimes \Gamma_{K|F}$ for any finite factor group $\Gamma_{K|F}=W_F/W_K,$ such that the normal subgroup $W_K$ acts trivially on $\wh G.$ The reason is that in all these cases the set $Z^1(W_F,\wh G)$ of 1-cocycles $\vf=\vf_2$ is in bijection with $Hom_{W_F}(W_F,\; \wh G\rtimes \Gamma_{K|F})$ via
$$  \vf\in Z^1(W_F,\wh G) \longmapsto \{w\mapsto \vf(w)\cdot \overline w \in \wh G\rtimes \Gamma_{K|F}\}.$$
In particular $Z^1(W_F,\wh G) = Hom(W_F,\wh G)$ if $W_F$ itself acts trivially on $\wh G.$
\enddefinition
If $\phi'\sim \phi$ are two equivalent (in the sense of 3.1(ii)) L-homomorphisms then the set $\{x\in\wh G\;|\; \;^x\phi=\phi'\}$ is always a coset $x_0C_{\wh G}(Im(\phi)).$
An important observation is

\proclaim{3.2 Proposition}Let $\phi=(\phi_1,\phi_2)$ be an L-homomorphism. Then the $\wh G$-centralizer
$C_{\wh G}(Im(\phi))$ is a subgroup of $\wh G$ which is reductive again but need not be connected.
\endproclaim

\demo{Proof} (Without proof the result is mentioned in {\bf [A],p.200}.) Let $\phi_1,\;\phi_2,\;\varphi_2$ be as explained in the definition. Then we have
$$Im \phi_1\subset \wh H_2:=C_{\wh G}(Im \phi_2)=\{x\in\wh G;\; x\varphi_2(w)\;^wx^{-1}
=\varphi_2(w)\;\text{for all}\; w\},$$
and actually $Im\phi_1\subseteq \wh H_2^0$ because $\phi_1$ is rational and $SL_2(\C)$ is connected. Furthermore:
$$ C_{\wh G}(Im \phi) = C_{\wh G}(Im \phi_1)\cap C_{\wh G}(Im \phi_2) =C_{\wh H_2}(Im \phi_1),$$
hence for the subgroups of 1-components we obtain:
$$C_{\wh G}(Im \phi)^0 = C_{\wh H_2}(Im\phi_1)^0=C_{\wh H_2^0}(Im \phi_1)^0.\tag{\bf 2}$$
As to the second equality the inclusion $\supseteq$ is obvious and the opposite inclusion follows from $C_{\wh H_2}(Im\phi_1)^0\subseteq \wh H_2^0.$
Since a linear group is reductive precisely when its 1-component is
reductive, we may restrict to the 1-components. In {\bf [K]10.1.1} and
{\bf [H]5.2} resp. it is proved that $\wh H_2=C_{\wh G}(Im\;\phi_2)$ is a
reductive group. The map $\phi_1:SL_2(\C)\rightarrow \wh
H_2^0\subset \wh G$ is rational, hence $Im\,\phi_1\subset H_2^0$ is a reductive subgroup of the reductive group $H_2^0.$ But in  {\bf [BV]}, proof of Proposition 2.4, it is mentioned that the centralizer of a reductive subgroup of a reductive group is again reductive. Therefore everything follows from {\bf (2)}.

Alternatively we remark that obviously
$$   \phi_1: SL_2(\C)\rightarrow D\wh H_2^0,\qquad d\phi_1:sl_2(\C)\rightarrow Lie(D\wh H_2^0),$$
where $D\wh H_2^0$ is the derived subgroup. Then from Prop.2.4 of loc.cit. we see directly that
$$\wh H_3:=\{x\in D\wh H_2^0\;|\; Ad(x)\circ Y=Y,\;\forall\;Y\in Im(d\phi_1)\}$$
is a reductive subgroup of $D\wh H_2^0.$ But $SL_2(\C)$ is simply connected and therefore from {\bf Bourbaki, LIE III,\S6.1, Thm 1}
we conclude that also $\wh H_3 = C_{D\wh H_2^0}(Im(\phi_1))\subset C_{\wh H_2^0}(Im(\phi_1)).$ Finally $\wh H_2^0$ modulo $D\wh H_2^0$ is generated by the center, hence together with $\wh H_3$ also the larger group $C_{\wh H_2^0}(Im(\phi_1))$ is reductive.\qed\enddemo


\bigskip
The background for all further considerations is:

\proclaim{3.3 Local Langlands conjecture}

\noindent Let $\Phi(G)$ be the set of $\wh G$-conjugacy classes of L-homomorphisms
$$   \phi: SL_2(\C)\times W_F\rightarrow \;^LG$$
which are our {\bf L-parameters.} Then there is a well defined bijection
$$   \Phi(G)\ni [\phi] \mapsto \Pi_{[\phi]} \tag{1}$$
between $\Phi(G)$ and the set of L-packets $\Pi_{[\phi]}$ (for the size of $\Pi_{[\phi]}$ see section{\bf 7.}) which form a partition of $Irr(G),$ the set of irreducible admissible representations of $G$ modulo equivalence. In other words the packets $\Pi_{[\phi]}$ are obtained as the fibers of a well defined surjective reciprocity map
$$  \pi\in Irr(G) \longmapsto [\phi](\pi)\in \Phi(G).\tag{2}$$

\noindent The subset $Irr(G)_{temp}$ of {\bf tempered} representations is precisely the union of the tempered L-packets corresponding (conjecturally) to the parameters $[\phi]$
such that the 1-cocycle
$$ \vf_2: W_F\rightarrow \wh G\qquad   \in Z^1(W_F,\wh G)$$
has relatively compact image in $\wh G$, i.e. the image has compact closure.\endproclaim
We recall that $\phi=(\phi_1,\phi_2)$ and the condition of temperedness depends only on the cocycle $\vf_2$
defining $\phi_2.$
Also note that the action of $W_F$ on $\wh G$ is via a finite quotient $W_F/W_K=Gal(K|F).$
And the cocycle condition is: $\vf_2(w_1w_2)=\vf_2(w_1)\;^{w_1}\vf_2(w_2).$
Therefore $\vf_2(W_Kw)=\vf_2(W_K)\vf_2(w),$
and the image of $\vf_2$ consists of finitely many cosets with respect to the subgroup $\vf_2(W_K).$
For temperedness it is therefore necessary and sufficient that this subgroup of $\wh G$ is relatively compact.
The conjecture and {\bf 1.3} suggest:
\definition{3.4 Definition}

{\bf (i)} The L-parameter $[\phi]\in \Phi(G)$ is called {\bf tempered} if for $\phi\in[\phi]$ the image of the 1-cocycle $\vf_2$ is relatively compact.

{\bf (ii)} An {\bf L-parameter standard-triple} $(P,[\tphi]_M,\nu)$ for $G$ consists of a standard $F$-parabolic subgroup, a tempered L-parameter $[\tphi]_M\in \Phi(M)$ of the standard Levi subgroup $M=M_P$, and a real parameter $\nu\in \aa_M^*$ which is regular with respect to $P.$
\enddefinition

\proclaim{3.5 Langlands classification of L-parameters}

\noindent There is a well defined bijection {\bf (for a precise definition see 4.6 below)}
$$   \{(P,[\tphi]_M,\nu)\}  \leftrightarrow  \Phi(G)      $$
which realizes the L-parameters for $G$ in terms of L-parameter standard-triples.\endproclaim

\definition{Remark} Note that we have here the same data $(P,\nu)$ as in the Langlands classification of irreducible representations in 1.4. This offers a possibility to extend a tempered reciprocity map
$ \sigma_M\in Irr(M)_{temp}\mapsto [\tphi]_M(\sigma_M)\in \Phi(M)_{temp}$ for the standard Levi groups to a larger map:
$$  (P,\sigma_M,\nu)\in Irr(G) \mapsto (P,[\tphi]_M(\sigma_M),\nu) \in \Phi(G),$$
which by its very definition takes nontempered representations to nontempered L-parameters.\enddefinition
\bigskip
{\bf 4. Twisting L-homomorphisms and a precise statement of 3.5}

For our reductive group $\bG|F$ we proceed with the subgroups $(\bP_0,\bM_0)$ and $(\bB,\bT)$ as we have fixed them in {\bf 2.3.2}. First of all we will give another characterization of the bijection $\bM\leftrightarrow\wh M$ between semi-standard Levi subgroups.

\proclaim{4.1 Proposition}
The bijection {\bf 2.3 (5a)} is an identification
$$   \bT\subset\bM\leftrightarrow \wh M\supset\wh T$$
between semi-standard Levi subgroups, which is characterized by the
equality
$$  X_*(Z(\wh M)^0) = X^*(\bM)\tag{1}$$
of sublattices in $X_*(\wh T)= X^*(\bT).$
\endproclaim
Before going to prove this we note that $X^*(\bM) \subseteq X^*(\bT)$ is not an inclusion in the proper sense but we have:\footnote {For a more general result see [KS],Lemma 13. We thank P.Schneider for referring us to that result.}
\proclaim{4.1.1 Lemma} The restriction from $\bM$ to $\bT$ induces an injection
$$  X^*(\bM) \cong X^*(\bT)^{W(\bT,\bM)}\subseteq X^*(\bT).$$
\endproclaim

\demo{Proof} We consider $\bM=\bG$ any connected reductive group with maximal torus $\bT$ and $W= N_\bG(\bT)/\bT.$ The restriction of characters (in the sense of algebraic groups) from $\bG$ to $\bT$ is injective because the natural map $\bT\rightarrow \bG/ D\bG$ is surjective (where $D\bG$ denotes the derived subgroup), and the image consists of $W$-invariant characters since $W$ acts via conjugation. So we have
$$  X^*(\bG)\hookrightarrow   X^*(\bT)^W.$$
Conversely let $\chi\in X^*(\bT)^W$ and let $\alpha^\vee\in X_*(\bT)$ be any coroot. Then $s_\alpha(\alpha^\vee) = -\alpha^\vee$ and the $W$-invariance of $\chi$ imply:
$$  \chi(\alpha^\vee) = \chi(-\alpha^\vee),\qquad (\chi\circ\alpha^\vee)^2 \equiv 1.$$
But we have $\chi\circ\alpha^\vee\in Hom(GL_1,GL_1)\cong \Z$ without torsion, hence $\chi\circ\alpha^\vee\equiv 1.$ So we see that $\chi\in (Q^\vee)^\perp$ where $Q^\vee\subset X_*(\bT)$ is the sublattice which is generated by the coroots. Moreover it is well known that

$\bT\cap D\bG$ is the maximal torus in the derived subgroup,

$Q^\vee \subseteq X_*(\bT\cap D\bG)$ is of finite index.

\noindent But this implies
$$  \overline F^\times\otimes Q^\vee =\overline F^\times \otimes X_*(\bT\cap D\bG) \cong \bT\cap D\bG,$$
where the isomorphism is due to {\bf [DM],0.20}. Therefore $\chi\in (Q^\vee)^\perp$ implies that $\chi$ is trivial on $\bT\cap D\bG,$ and from $\bT/\bT\cap D\bG \cong \bG/D\bG$ we see now that $\chi$ is extendable onto $\bG.$\qed
\enddemo
\demo{Proof of {\bf 4.1}} The semi-standard Levi group $\wh M$  is recovered as
$$   \wh M = C_{\wh G}((\wh T^{W(\wh T,\wh M)})^0)\quad\text{hence}\quad Z(\wh M)^0=(\wh T^{W(\wh T,\wh M)})^0,$$ because $Z(\wh M)^0\subseteq (\wh T^{W(\wh T,\wh M)})^0$ is obvious,
and this implies
$$  X_*(Z(\wh M)^0)= X_*(\wh T^{W(\wh T,\wh M)})= X_*(\wh T)^{W(\wh T,\wh M)}\subseteq X_*(\wh T).$$
But by construction we have $X_*(\wh T)=X^*(\bT)$ including the equality of subsets $\Delta_*(\wh T,\wh B)=\Delta^*(\bT,\bB),$ and dually we have the equality $\Delta^*(\wh T,\wh B)=\Delta_*(\bT,\bB).$ So we see that the roots / coroots give rise to the same reflections of $X_*(\wh T)=X^*(\bT)$ and therefore we obtain the same Weyl groups $W(\wh T,\wh G)=W(\bT,\bG)\subset Aut(X_*(\wh T))=Aut(X^*(\bT)).$ Finally under this identification we also have

\noindent $W(\wh T,\wh M)= W(\bT,\bM)$ as we see from the remarks following {\bf 2.3 (5a)}. So we arrive at
$$  X_*(Z(\wh M)^0)= X_*(\wh T)^{W(\wh T,\wh M)}= X^*(\bT)^{W(\bT,\bM)}\cong X^*(\bM)$$
where the last isomorphism is due to {\bf 4.1.1}.

\noindent {\it Finally we remark that for a given semi-standard Levi group $\bM\supset\bT$ there can exist only one Levi group $\wh M\supset\wh T$ such that $X_*(Z(\wh M)^0)\subseteq X_*(\wh T)$ identifies with the image $X^*(\bT)^{W(\bT,\bM)}$ of $X^*(\bM)$ in $X^*(\bT)=X_*(\wh T).$} For this it is enough to see that $\wh M$ can be recovered from $X_*(Z(\wh M)^0) \subseteq X_*(\wh T).$ But the natural identification $\C^\times\otimes X_*(\wh T)\cong \wh T$ will take $\C^\times\otimes X_*(Z(\wh M)^0)$ to the subtorus $Z(\wh M)^0\subseteq\wh T,$ and $\wh M$ is recovered as $\wh M= C_{\wh G}(Z(\wh M)^0).$ So we see that our identification {\bf (1)} determines a map $\bM\mapsto \wh M$ which must be the bijection {\bf (5a)} in 2.3.
\qed\enddemo

\proclaim{4.2 Corollary} From {\bf (1)} we obtain natural isomorphisms

$$ \CD
Z(\wh M)^0@>\sim>> \C^\times\otimes X_*(Z(\wh M)^0)@>=>> \C^\times\otimes X^*(\bM)\\
@VV\text{incl} V             @VVV                                @VV\text{res} V\\
\wh T     @>\sim>> \C^\times\otimes X_*(\wh T)     @>=>> \C^\times\otimes X^*(\bT),
\endCD $$
where the isomorphisms from right to left are: $\lambda^{\otimes \chi^\vee}\mapsto \chi^\vee(\lambda).$ Moreover the map in the lower line is
equivariant with respect to the action of Weyl groups: $N_{\wh G}(\wh T)/\wh T \cong W_*(\wh T,\wh G)=W^*(\bT,\bG)$ which means that $ \lambda^{\otimes w\circ\chi^\vee}\mapsto \wh n_{\wh w}\cdot\chi^\vee(\lambda)\cdot \wh n_{\wh w}^{-1}.$\qed

\endproclaim

Since the bijection $\bM\leftrightarrow \wh M$ of semi-standard Levi-subgroups is $\mu(\Gamma)$-$\Gamma$ equivariant, we conclude from {\bf (1)} that
$$  X_*(Z(\;^LM)^0) = X_*(Z(\wh M)^0)^\Gamma = X^*(\bM)^{\mu(\Gamma)}\tag{2}$$
if $\bM\leftrightarrow \wh M$ are $\mu(\Gamma)$-$\Gamma$-stable groups, and $Z(\;^LM)^0$ is the connected center of $\wh M\rtimes W_F.$  Next we look at semi-standard $F$-Levi groups.
\proclaim{4.3 Corollary} Let $\bM|F \subseteq\bG|F$ be a semi-standard $F$-Levi-subgroup of $\bG$. Then from {\bf (1)} we obtain the identification
$$  X_*(Z(\;^LM)^0) =X^*(\bM)_F,\quad \chi^\vee \leftrightarrow \chi \tag{i}$$
of $\Gamma$-invariant sublattices, which, as in {\bf 4.2}, turns into an identification
$$Z(\;^LM)^0\cong \C^\times\otimes_\Z X_*(Z(\;^LM)^0)=\C^\times\otimes_\Z
X^*(\bM)_F, \quad \chi^\vee(\lambda)\leftrightarrow \lambda^{\otimes \chi}.\tag{ii}$$
Using the minimal $F$-Levi group $\bM_0,$ this embeds into
$$Z(\;^LM_0)^0\cong \C^\times\otimes_\Z X_*(Z(\;^LM_0)^0)=\C^\times\otimes_\Z
X^*(\bM_0)_F.\tag{iii}$$
which is now acted on by $\;_F\wh W\cong \;_FW$ in an equivariant way.
\endproclaim
\demo{Proof} We use {\bf (2)} and the Corollary in 2.2.3 which tells us that on $X^*(\bM)\hookrightarrow X^*(\bT)$ the $\mu(\Gamma)$-action and the usual $\Gamma$-action agree if $\bM$ is an $F$-group.

\noindent The last statement follows from the corresponding statement in {\bf 4.2} together with {\bf 2.2.7} and {\bf 2.3.4 (ii)}.
\qed\enddemo

\noindent Because $Z(\;^LM)^0$ is a complex torus we have the well defined subgroup of hyperbolic elements, and the identification {\bf 4.3(ii)} restricts to
$$ Z(\;^LM)^0_{hyp} \cong (\R_+)^\times\otimes X^*(\bM)_F \cong \aa_M^*,\quad \chi^\vee(q^\beta)\leftrightarrow q^{\beta\otimes\chi}\leftrightarrow \beta\otimes\chi, \tag{3}$$
where the second isomorphism has been considered in {\bf 1.2 (2)}, and again {\bf (3)} is subordinate to the  $\;_F\wh W$-$\;_FW$ equivariant identification
$$ Z(\;^LM_0)^0_{hyp} \cong (\R_+)^\times\otimes X^*(\bM_0)_F \cong \aa_{M_0}^*. \tag{4}$$
\definition{4.4 Definition} We will write the maps (3) / (4) as
$$  \nu\in \aa_M^* \mapsto z(\nu)\in Z(\,^LM)^0_{hyp},\qquad z(\beta\otimes\chi) = \chi^\vee(q^\beta).$$
For the Langlands classification of L-parameters they play the same role as $\nu\in\aa_M^*\mapsto \chi_\nu$ (considered in {\bf 1.2}) for the Langlands classification of representations.
\enddefinition
Just as we may twist an irreducible representation $\pi$ by an unramified character $\chi$ we have now the possibility to twist an L-parameter $[\phi]\in \Phi(G)$ by an element $z\in Z(\,^LG)^0$ of the central torus. Using the exponent $w\mapsto d(w)$ of {\bf 3.1 (iii)} we have:



\definition{4.5 Definition (see [A],p.201)} If
$$  \phi=(\phi_1,\phi_2):SL_2(\C)\times W_F\rightarrow \;^LG,\qquad z\in Z(\;^LG)^0= (Z(\wh G)^\Gamma)^0$$
is a pair consisting of an L-homomorphism and an element in the connected center of the L-group,
then the {\bf twisted L-homomorphism} is:
$$  \phi_z(A,w):= \phi(A,w)\cdot z^{d(w)}= \phi_1(A)\phi_2(w)\cdot z^{d(w)}\tag{7}$$
and the {\bf twisted L-parameter} is:
$$ [\phi]_z := [\phi_z] \in \Phi(G).$$\enddefinition
This means that $\phi_1$ is left unchanged, but the cocycle $\vf_2$ is replaced by the twisted cocycle
$$  \{w\mapsto \vf_2'(w):= \vf_2(w)\cdot z^{d(w)}\} \in Z^1(W_F,\wh G).$$

Using 4.4 and 4.5 the Langlands classification 3.5 of L-parameters has the more precise formulation:

\proclaim{4.6 Langlands classification of L-parameters revisited} There is a bijection
$$  (P,[\tphi]_M,\nu)\mapsto [\phi]_G\in\Phi(G)\tag{9}$$
between L-parameter standard-triples and L-parameters, where

{\bf (i)} $\phi = (\tphi)_{z(\nu)}$ is the twisted L-homorphism (taken on the level of $M$) as it is defined in {\bf (7)}, and $[\phi]_G$ is obtained from $\phi$ via $\;^LM\subseteq \;^LG,$

{\bf (ii)} $z(\nu)\in Z(\;^LM)^0_{hyp}$ is assigned to $\nu\in\aa_M^*.$

{\bf (iii)} In terms of {\bf 5.3 (i)} below, the element $z(\nu)$ is recovered from $[\phi]_G$ as $z(\nu)=z(\phi)$ for an appropriate representative $\phi\in [\phi]_G.$
\endproclaim

Parallel to {\bf 1.4*} we have
\proclaim{4.6* A second realization} A different bijection between L-parameter standard-triples and L-parameters
can be realized as follows:
$$  (P,[\tphi]_M,\nu)\mapsto [\phi]_G\in\Phi(G)\tag{9*}$$
where $\phi = (\tphi)_{z(-\nu)}$ is the twist of $\tphi$ by $z(-\nu)=z(\nu)^{-1}\in Z(\,^LM)^0_{hyp}$ and $[\phi]_G$ is obtained from $\phi$ via $\;^LM\subseteq \;^LG.$\endproclaim

\noindent {\bf Remark:} In terms of {\bf 5.3 (i)} below, the element $z(\nu)$ is recovered now from $[\phi]_G$ as $z(\nu)=z(\phi(w_{-1}))$ for an appropriate representative $\phi\in [\phi]_G,$
where $w_{-1}\in W_F$ is a geometric Frobenius-lift. The proof of 4.6* uses the obvious reformulation of {\bf 5.3} where Frobenius lifts $w_1$ are replaced by geometric Frobenius lifts $w_{-1}.$

As a background information (which plays no role in the following) we mention another desideratum to be expected of the Langlands correspondence:

\proclaim{4.7 Conjecture on unramified twists} Consider the isogeny $z\in Z(\;^LG)^0\mapsto \chi_z\in X_{ur}(G)$ as it is obtained by combining {\bf 4.3 (ii)} and {\bf 1.2 (1)}. Then
the conjectural Langlands map $[\phi]\in\Phi(G)\mapsto \Pi_\phi\subset Irr(G)$ of 3.3, is expected to have the {\bf twist property}
$$ \Pi_{\phi_z} = \Pi_{\phi}\otimes \chi_z,\tag{i}$$
which is an identity of L-packets. In particular for $z=\chi^\vee(q^s)\leftrightarrow q^{s\otimes\chi}\in\C^\times\otimes X^*(\bG)_F$ we should have
$$  \Pi_{\phi\cdot(\chi^\vee\circ\pa .\pa^s)} = \Pi_\phi\otimes(| .|_F^s\circ\chi),\tag{ii}$$
if $\chi\in X^*(\bG)_F$ is a rational character, $\chi^\vee\in X_*(Z(\,^LG)^0)$ the corresponding cocharacter and $s\in\C.$

\noindent In terms of the reciprocity map {\bf 3.3(2)} the property {\bf (i)} means that
$$  [\phi](\pi\otimes\chi_z) = ([\phi](\pi))_z  \in \Phi(G).$$
\endproclaim
\noindent We note that {\bf (i)} implies {\bf (ii)} because on one hand $z=s\otimes\chi\in \C\otimes X^*(\bG)_F$ is mapped to $\chi_z=|\chi|_F^s$ as we see from {\bf 1.2 (1)}, and on the other hand substituting $z=\chi^\vee(q^s)$ into (7) we obtain:

$$ z^{d(w)}=\chi^\vee(q^s)^{d(w)} = \chi^\vee((q^{d(w)})^s)= \chi^\vee\circ\pa w\pa^s.$$

A more general conjecture of Langlands is stated in {\bf [Bo1], 10.2 and 10.3(2)}:
We have a well defined action
$$    H^1_c(W_F, Z(\wh G)) \times  \Phi(G)\rightarrow \Phi(G) ,\qquad [\alpha]\cdot [\phi]:= [\alpha\cdot\phi]$$
and a homomorphism
$$   [\alpha]\in H^1_c(W_F, Z(\wh G)) \mapsto \chi_\alpha\in X_c(G),$$
and the conjectural Langlands map is expected to have the {\bf twist property}
$$  \Pi_{\alpha\cdot\phi} = \Pi_\phi\otimes \chi_\alpha, \qquad [\phi](\pi\otimes\chi_\alpha)=[\alpha]\cdot ([\phi](\pi)).  \tag{8}$$
It is in fact possible to show that {\bf 4.7} is subordinate to {\bf (8)}. (see {\bf [Lvor], justification of 2.3.2}). This justifies the definition {\bf (7)} of twisting L-homomorphisms by elements $z\in Z(\;^LG)^0$ from the connected center of the L-group.

\noindent {\bf Remark:} In {\bf 4.7} we have obtained a natural isogeny
$Z(\;^LG)^0 = \C^\times\otimes X^*(\bG)_F\twoheadrightarrow X_{ur}(\bG(F)).$ Using the maximal $F$-split torus $\bS_G$ which comes as a quotient of $\bG$ one may recover this map as
$$ \CD
\wh S_G @> \sim >> X_{ur}(\bS_G(F))\\
@VV\cong V               @VV\text{inf} V\\
Z(\,^LG)^0                 @>>> X_{ur}(\bG(F))
\endCD $$
where the vertical maps are natural isomorphism / inflation resp. and where the upper row is the unramified part of the Langlands correspondence for tori. We could also use $\bA_G|F,$ the maximal $F$-split torus in the center of $\bG,$ and then by a projection /restriction map go down to the line $\wh A_G\cong X_{ur}(\bA_G(F)).$

So far we have obtained a procedure how to produce L-parameters of $G$ out of L-parameter standard-triples.
Now we are going to construct the converse map of {\bf (9)}.
\bigskip

{\bf 5. Two basic invariants associated to an L-homomorphism}

\noindent {\bf 5.1 Polar decomposition in $\;^LG.$}

We recall that any semisimple element $s$ in a complex reductive
group has a canonical {\bf polar decomposition} $s=s_e\cdot s_h$ into
commuting elliptic and hyperbolic factors resp., such that for all
rational characters $\chi$ of diagonalizable subgroups containing $s$ the
values $\chi(s_e),\chi(s_h)\in \C^\times$ will be on the unit circle
and in $\R_{>0}$ resp. Since a diagonalizable group is the direct product
of its connected component (torus) and a finite group, {\bf the hyperbolic
factor $s_h$ will always be an element of the connected component.}

{\bf Note} here that a semisimple element of an algebraic group can
be always embedded into a closed subgroup which is diagonalizable.
{\bf [Bo2]8.4}.

\noindent Now we consider the group $\;^LG=\wh G\rtimes W_F.$ For
every finite Galois extension $K|F$ which splits $\bG$ we can form
the algebraic factor group $\;^LG/W_K = \wh G\rtimes \Gamma_{K|F},$
and we will simply speak here of algebraic factor groups $\wh
G\rtimes\Gamma_*$ of $\;^LG.$ An element $s\in\;^LG$ is called
semisimple, hyperbolic, elliptic, resp. if the projection of $s$ to
all algebraic factor groups of $\;^LG$ has this property.  Now if
$s\in \;^LG$ is semisimple we have its polar decomposition
$$ s_*=s_{*,h}\cdot s_{*,e} $$
in each algebraic factor group $\wh G\rtimes\Gamma_*.$ But
$s_{*,h}=s_h\in\wh G$ does not depend on $\Gamma_*.$

\noindent (We may embed $s_*$ into a diagonalizable subgroup
$D\subset \wh G\rtimes\Gamma_*$ and then $s_{*,h}$ is in the
connected component of $D$ hence in $\wh G.$)

\noindent Therefore the commutator $[s,s_h]$ is in
$Ker\{\;^LG\rightarrow \wh G\rtimes\Gamma_*\}$ for all $\Gamma_*,$
hence $s$ and $s_h$ commute and $s=s_h\cdot\e,$ where $\e=
s_h^{-1}s,$ is the polar decomposition of $s.$ In particular we see
that the polar decomposition of $s$ is fixed already by looking at a
single group $\wh G\rtimes\Gamma_*$ and that a semisimple hyperbolic
element of $\;^LG$ is the same as a semisimple hyperbolic element of
$\wh G.$

\noindent We want to apply this in the situation of {\bf 3.1} to the
semisimple elements $\phi_2(w)=(\varphi_2(w),w)\in \;^LG.$ Then we
obtain a polar decomposition
$$s:=(\vf_2(w),w)=\phi_2(w) = s_h\cdot s_e$$
with hyperbolic part $s_h\in \wh G$ and elliptic part
$s_e=(s_h^{-1}\varphi_2(w),\;w)\in \;^LG.$ Since $\phi_2(w)$ and
$s_h$ must commute, we obtain
$$  s_h\cdot \varphi_2(w)= \varphi_2(w)\cdot w(s_h)\quad\in\wh G.$$
So we can neither conclude that $\varphi_2(w)$ is semisimple nor
that $s_h$ is the hyperbolic part of $\varphi_2(w)\in\wh G.$

\noindent {\bf But if $s=(\vf_2(w),w)\in\;^LG$ is such that $\;^w\vf_2(w) =\vf_2(w)\in \wh G,$} then $s$ is semisimple if and only if $\vf_2(w)$ is, and the polar decomposition is then:
$$  s_h =\vf_2(w)_h,\qquad s_e= (\vf_2(w)_e,w).$$

We can now reformulate the temperedness criterion {\bf 3.4(i)} as
follows:

\proclaim{5.2 Lemma} For a Langlands parameter $\phi:SL_2(\C)\times
W_F\rightarrow \;^LG$ the following is equivalent:

{\bf (i)} $\phi$ is tempered, which means that
$\varphi_2(W_F)\subset \wh G$ has compact closure,

{\bf (ii)} If $K|F$ is a finite extension such that $W_K\subset W_F$ acts trivially on $\wh G,$
then $\vf_2(W_K)$ consists entirely of elliptic semisimple elements,

{\bf (iii)} The semisimple elements
$\phi_2(W_F)\subset \;^LG$ are all elliptic.
\endproclaim

\demo{Proof} Assume that $Im(\phi_2)$ consists entirely of elliptic
elements. We consider a subgroup $W_K\subset W_F$ of finite index
such that $W_K$ acts trivially on $\wh G.$ Then the components of
$\phi_2(w)=(\varphi_2(w),w)$ commute if $w\in W_K,$ and therefore
together with $\phi_2(w)$ also $\varphi_2(w)$ is semisimple
elliptic. Moreover the cocycle $\varphi_2$ restricted to $W_K$ is a
homomorphism, hence $\varphi_2(W_K)\subset \wh G$ is a subgroup
which entirely consists of semisimple elliptic elements. Now let
$I_K\subset W_K$ be the inertia subgroup. It is compact and totally
disconnected whereas in $\wh G$ we find a neighborhood  of 1 which
besides $\{1\}$ does not contain any subgroup. Therefore
$\varphi_2(I_K)\subset \wh G$ is a finite subgroup. Consider
$s_1:=\varphi_2(w_1)$ where $w_1\in W_K$ is a Frobenius-lift. Then
$s_1$ normalizes the finite group $\varphi_2(I_K),$ hence a certain
power $s_1^d$ will centralize this group. Consider now
$W_L=I_K\rtimes <w_1^d>\quad \subset W_K.$ By construction we have
$$  \varphi_2(W_L) = \varphi_2(I_K)\times <s_1>,$$
the direct product of a finite group and a cyclic group which is
generated by a semisimple elliptic element. Therefore
$\varphi_2(W_L)$ is relatively compact and since it is of finite
index in $\varphi_2(W_K)$ the same holds for $\varphi_2(W_K).$ But
the image $\varphi_2(W_F)$ consists of finitely many
$\varphi_2(W_K)$-cosets, so it is relatively compact too.

\noindent {\bf Conversely} assume now that $\varphi_2(W_F)\subset\wh
G$ has compact closure. The same is then true for $\varphi_2(W_K).$
Similar as above we see that $\varphi_2(W_K)\subset\wh G$ is a group
which consists entirely of semisimple elements. Therefore compact
closure implies that all these elements are elliptic. Now we
consider $\phi_2(W_K)=\{(\varphi_2(w),w)\;|\;w\in W_K\}.$ Because
here $\varphi_2(w)$ and $w$ always commute, we see that
$\phi_2(W_K)$ consists entirely of elliptic semisimple elements of
$\;^LG.$ Since all elements of $\phi_2(W_F)$ are semisimple and
$\phi_2(W_K)$ is a subgroup of finite index we conclude that all
elements of $\phi_2(W_F)$ are elliptic. \qed\enddemo

In view of the Lemma the following Proposition will be the key for
the classification result.
\proclaim{5.3 Proposition}  Let $\phi:SL_2(\C)\times
W_F\rightarrow \;^LG$ be any $L$-homomorphism, and let $w_1\in W_F$ be a
Frobenius lift which means the exponent is $d(w_1)=1.$ Then:

{\bf (i)} The hyperbolic part $z$ of the semisimple element $s'=
\phi_2(w_1)=(\varphi_2(w_1),w_1)\in\;^LG$ is contained in the
centralizer $C_{\wh G}(Im(\phi))^0$ and does not depend on the
choice of $w_1.$ {\bf Actually $z=z(\phi)$ is even contained in $Z(C_{\wh
G}(Im(\phi))^0$ which is the maximal torus in the center of the
reductive group $C_{\wh G}(Im(\phi))$.}

{\bf (ii)} The element $z$ from (i) is the uniquely determined
semisimple hyperbolic element $z\in C_{\wh G}(Im(\phi))^0$ such that
$$   \phi_t(A,w):= \phi(A,w)\cdot z^{-d(w)}$$ is a tempered
$L$-homomorphism.

{\bf (iii)}  The L-homomorphism $\phi_t=(\phi_{t,1},\phi_{t,2})$ has
the property that $\phi_{t,2}(w)=\phi_2(w)_e$ is always the elliptic
part of the semisimple element $\phi_2(w)=(\varphi_2(w),w)\in
\;^LG.$
\endproclaim The Proposition is mentioned in {\bf [A],p.201} in a
less precise form. The main result is (i), which is a variation of
{\bf [H], 5.1.} We follow the proof which is given there.

\demo{Proof} We begin by explaining that {\bf (i) implies (ii)}.
First we show that $\phi_t$ is indeed a tempered L-homomorphism if
$z$ has the properties announced in (i). According to the Lemma we
need to show that $\phi_{t,2}(W_F)$ consists entirely of elliptic
elements. Here we use the following

\proclaim{Claim} If $\ve\in\;^LG$ is a semisimple elliptic element
which normalizes $\phi_2(I_F)\subset \;^LG$ then the whole coset
$\ve\phi_2(I_F)$ consists of semisimple elliptic
elements.\endproclaim \demo{Proof} We take an algebraic quotient
$\;^LG_*=\wh G\rtimes\Gamma_*$ corresponding to a finite Galois
extension $K|F$ which splits $\bG.$ Then $\phi_2$ induces
$\phi_{2,*}:W_F\rightarrow \;^LG_*.$ As in the proof of the previous
Lemma we see that $\phi_{2,*}(I_K)=\varphi_2(I_K)$ is finite, hence
$\phi_{2,*}(I_F)$ is finite too, and is normalized by the projected
element $\ve_*\in\;^LG_*.$ We find then a power $\ve_*^d$ which
centralizes $\phi_{2,*}(I_F)$ and therefore
$\phi_{2,*}(I_F)<\ve_*^d>$ consists entirely of semisimple elliptic
elements. Since $\phi_{2,*}(I_F)<\ve_*>\quad\supset
\phi_{2,*}(I_F)<\ve_*^d>$ is of finite index, the same is true for
$\phi_{2,*}(I_F)\cdot <\ve_*>$ hence also for $\phi_2(I_F)\ve$
because we have worked with an arbitrary quotient $\;^LG_*.$
\qed\enddemo

Now assuming {\bf (i)} we have
$$  \phi_2(w_1) = z\cdot s_e\quad \in \;^LG$$
where $z$ is hyperbolic and $s_e\in \;^LG$ is semisimple elliptic.
Since $z$ centralizes $Im(\phi_2)$ we see that together with
$\phi_2(w_1)$ also $s_e$ will normalize $\phi_2(I_F).$ Moreover by
definition of $\phi_t$ we get:
$$\phi_{t,2}(w)=\phi_2(w)\cdot z^{-d(w)}=\phi_2(ww_1^{-d(w)})\cdot
s_e^{d(w)},$$ where $ww_1^{-d(w)}\in I_F.$ Therefore using the Claim
we see that $\phi_{t,2}(w)$ is always elliptic.

\noindent {\bf Conversely} if $\phi_t$ is tempered then in
particular $\phi_{t,2}(w_1) = \phi_2(w_1)\cdot z^{-1}$ must be
elliptic. So if we already know that $z\in C_{\wh G}(Im \phi)^0$ and
that $z$ is hyperbolic, then $z$ must be the hyperbolic factor of
the semisimple element $\phi_2(w_1).$ {\bf So we have seen that (i) implies (ii).}

\noindent {\bf As to (iii)} the argument is the same as before since $\phi_t(I_2,w)$ must be semisimple elliptic.

\noindent {\bf (i)} Let $\phi_*$ be the combined map
$$  SL_2(\C)\times W_F\rightarrow \;^LG\rightarrow \wh G\rtimes \Gamma_*$$
and consider $s'=\phi_*(w_1).$
Our first aim is to show that the hyperbolic part $z$ of $s'$ is in the
centralizer $C_{\wh G}(Im(\phi_*)).$
Let $I_F\subset W_F$ be the inertia group. Because $I_F$ is compact and
totally disconnected we see that $\phi_*(I_F)\subset \wh G\rtimes\Gamma_*$
is a finite subgroup which is normalized by $s'.$ Therefore we find a power
${s'}^m$ which commutes with all elements from $\phi_*(I_F).$ But ${s'}^m\in\wh G\rtimes\Gamma_*$
and the group $\Gamma_*$ is finite, hence for an appropriate $m'$ we obtain
$(s')^{mm'}\in \wh G$ and $(s')^{mm'}$ commutes with $\phi_*(W_F)=\phi_*(I_F)\rtimes<s'>.$
Moreover $(s')^{mm'}\in \phi_*(W_F)$ commutes with $\phi_*(SL_2(\C))$, therefore
putting $l=mm'$ we have
$$  (s')^l\in C_{\wh G}(Im(\phi_*)).$$
If the order of $s'$ is finite then $\phi_*(W_F)$ is a finite group of semisimple elements and the hyperbolic part of $s'=\phi_*(w_1)$ is $z=1$ independently of the choice of $w_1.$

\noindent Now assume that $s'$ is of infinite order. The centralizer $\H:= C_{\wh G\rtimes\Gamma_*}((s')^l)$ is a reductive group (since according to {\bf [C]} Theorem 3.5.4 on p.93, the connected component $\H^0=C_{\wh G}((s')^l)^0$ is reductive) and $Im(\phi_*)\subset \H\subset \wh G\rtimes\Gamma_*.$ Because $(s')^l\in\H$ is a central semisimple element of infinite order, the maximal central torus $Z_{\H}$ of $\H$ is nontrivial and $(s')^{kl}\in Z_{\H}$ for a certain power $k.$ The torus $Z_{\H}$ being complex we find $s\in Z_{\H}$ such that $s^{kl}={s'}^{kl}.$ By construction
we see that $Z_{\H}$ and in particular the element $s$ commutes
with $Im(\phi_*)\ni s'.$ Since $Z_{\H}$ is connected we must have $Z_{\H}\subset C_{\wh G}(Im(\phi_*))^0.$ We deduce that
$$   s'= s's^{-1}\cdot s$$
is a decomposition into commuting semisimple factors where the first
one is of finite order. Therefore $s'$ and $s\in Z_{\H}$ have in their polar
decomposition the same hyperbolic factor $z.$ Together with $s$,
also $z\in Z_{\H}\subset C_{\wh G}(Im(\phi_*))^0.$ But obviously any $x\in C_{\wh G}(Im(\phi))$ commutes with $s'=\phi_*(w_1)$ and therefore it also commutes with the hyperbolic part $z$ of $s'.$ Hence $z$ is in the center of $C_{\wh G}(Im(\phi)),$ and even in the connected component of the center because it is hyperbolic.

\noindent {\bf Finally} we show that $z$ does not depend on the choice of $w_1.$ If we change $w_1$ then we obtain
$$   s'' = \iota\cdot s'\in \phi_*(W_F)\subset\wh G\rtimes\Gamma_*,$$
where $\iota\in \phi_*(I_F).$ Since $\phi_*(I_F)$ is a finite group which is normalized by $s'$, a certain power ${s'}^m$ will commute with all elements from $\phi_*(I_F).$  Therefore:
$$ {s''}^m=(\iota\cdot s')^m = \iota_1\cdot {s'}^m,\qquad {s''}^{mm'}=(\iota_1\cdot{s'}^m)^{m'}= \iota_1^{m'}\cdot {s'}^{mm'},$$
for all $m'$, and we find a power $l=mm'$ such that ${s''}^l={s'}^l.$ But $z$ commutes with $Im(\phi_*)$ and therefore
$$  (s''z^{-1})^l = {s''}^lz^{-l}= {s'}^lz^{-l}= (s'z^{-1})^l$$
is an elliptic semisimple element. Hence $s''z^{-1}$ is elliptic too, and $s''=(s''z^{-1})\cdot z$ is the polar decomposition of $s''.$ \qed\enddemo

Summarizing {\bf 5.3} we have a well defined map
$$  \phi \mapsto z(\phi)\in Z(C_{\wh G}(Im(\phi)))^0,\tag{1}$$
which associates a central hyperbolic element to a L-homomorphism $\phi.$
If we go to an equivalent L-homomorphism $x\phi x^{-1},$ where $x\in \wh G$, then of course $z(x\phi x^{-1})= xz(\phi)x^{-1}.$

\noindent In {\bf 5.3} we could have used also the second copy of $W_F$ in $SL_2(\C)\times W_F$ which comes from the Weil-Deligne group. (see {\bf 3.1(iii)}). This gives us another invariant
$$ z_*(\phi):= z(\phi(w_{1*})) = \phi_1\left(\bmatrix q^{1/2}&0\\ 0&q^{-1/2}\endbmatrix\right)\cdot z(\phi).\tag{1*}$$
As a product of two commuting hyperbolic elements $z_*(\phi)$ is again hyperbolic semisimple, so it is the hyperbolic part of $\phi(w_{1*}),$ and
$$ z_*(\phi)\in \phi_1(D(\C))\cdot Z(C_{\wh G}(Im(\phi)))^0,$$
where $D(\C)$ denotes the diagonal torus in $SL_2(\C).$

Therefore from {\bf 3.4, 5.2, 5.3} we see:
\proclaim{5.4 Corollary} $[\phi]\in \Phi(G)$ is a tempered L-parameter if and only if $z(\phi)=1,$
or equivalently $z_*(\phi)=\phi_1\left(\bmatrix q^{1/2}&0\\ 0&q^{-1/2}\endbmatrix\right).$
\qed\endproclaim

\noindent {\bf From 5.3 we want to recover}
$$   \phi = (\phi_{z^{-1}})_z,$$
but the last twist is only allowed if $z\in Z(\;^LG)^0$ (see above). Therefore
we search  for Levi-subgroups $L(\;^LG)\subseteq \;^LG=\wh G\rtimes W_F$  having the properties:
$$  Im(\phi)\subseteq L(\;^LG),\quad\text{and}\; z\in Z(L(\;^LG))^0.\tag{P1}$$
{\bf But this is equivalent to the property:}

\noindent {\bf (P1)*:} $\quad L(\;^LG)= C_{\;^LG}(\wh S)$ is the $\;^LG$-centralizer of a complex torus $\wh S$ such that $z\in\wh S\subset C_{\wh G}(Im(\phi))^0.$

\demo{Proof} To see the equivalence we may use the following weak form of {\bf [Bo1],3.5.}
\proclaim{Lemma} If $\wh S\subseteq \wh G$ is a torus such that the centralizer $C_{\;^LG}(\wh S)$ meets
every connected component of $\;^LG$, then $C_{\;^LG}(\wh S)= L(\;^LG)$ is a Levi subgroup, and conversely each Levi subgroup $L(\;^LG)$ is obtained in that way, because $L(\,^LG)$ is recovered from its connected center:
$$ L(\;^LG) = C_{\;^LG}( Z(L(\;^LG))^0).$$
\qed\endproclaim

Therefore, if $L(\;^LG)$ is a Levi subgroup such that (P1) is satisfied, then $\wh S:= Z(L(\;^LG))^0 = C_{\wh G}(L(\;^LG))^0$ is a torus which has the required properties.

\noindent {\bf Conversely} if $\wh S$ is a torus as in (P1)*, then $C_{\;^LG}(\wh S)\supseteq Im(\phi)$ and therefore it meets every connected component of $\;^LG.$ So we see that indeed
$C_{\;^LG}(\wh S)=L(\;^LG)$ is a Levi subgroup which satisfies (P1).\qed\enddemo

Obviously there is a unique minimal $\wh S$ such that $z\in\wh S\subset Z(C_{\wh G}(Im(\phi)))^0,$ namely
$\wh S = \overline{<z>}^0$ the connected component of the Zariski closure of the cyclic group $<z>.$ (This must contain $z$ because the element $z$ is hyperbolic.) Therefore:
\proclaim{5.5 Proposition} Assigned to an L-homomorphism $\phi$ we have not only a unique hyperbolic semisimple element $z=z(\phi)$ but also a unique {\bf maximal} Levi subgroup $L(\;^LG)_\phi$ such that the condition {\bf (P1)}holds, and these assignments are compatible with $\wh G$-conjugation. Thus we will have
$$  z(x\phi x^{-1})= x z(\phi) x^{-1},\qquad  L(\;^LG)_{x\phi x^{-1}}= x L(\;^LG)_\phi x^{-1},$$
for any $x\in \wh G.$
\endproclaim

\bigskip

{\bf 6. On the proof of 3.5 / 4.6}

It is not hard to see that the proof of 3.5 /4.6 will be finished if we can construct the converse map, which means
beginning from $[\phi]_G\in \Phi(G)$ we have to find the corresponding L-parameter standard triple. As a first step in that direction we have obtained already the Proposition 5.5. Now the construction will be completed by the following result:

\proclaim{6.1 Proposition} Let $[\phi]$ be a $\wh G$-conjugacy class of L-homomorphisms

\noindent $\phi:SL_2(\C)\times W_F\rightarrow \;^LG,$ and let $(z(\phi),\; L(\;^LG)_\phi)$ be the uniquely determined pair consisting of a semisimple hyperbolic element and a Levi subgroup of $\;^LG$ which is maximal with property (P1). Then:

{\bf (i)} It is possible to choose $\phi\in [\phi]$ such that moreover

\noindent {\bf (P2):} the group $L(\;^LG)_\phi$ is a {\bf relevant} standard Levi group in $\;^LG$, and the hyperbolic element $z(\phi)\in Z(L(\;^LG)_\phi)^0$ is regular with respect to the standard parabolic group $P(\;^LG)_\phi\supset L(\;^LG)_\phi.$

{\bf (ii)} Equivalently to {\bf (i)} we have
$$  L(\;^LG)_\phi =\;^LM,\qquad P(\;^LG)_\phi = \;^LP $$
with standard F-groups $(P,M)$ in $G=\bG(F),$  and under {\bf 4.3 (3)}
$$   z(\phi)\in Z(\;^LM)^0_{hyp} \leftrightarrow \nu(\phi)\in\aa_{M}^* $$
is regular with respect to $P.$

{\bf (iii)} If $\phi',\phi\in [\phi]$ are two conjugate L-homomorphisms satisfying (i), then we obtain identical triples
$$  \{z(\phi'),\;L(\;^LG)_{\phi'},\;P(\;^LG)_{\phi'}\} = \{z(\phi),\;L(\;^LG)_{\phi},\;P(\;^LG)_{\phi}\}$$
and $\phi'=x\phi x^{-1},$ where $x\in L(\wh G)_\phi =\wh{M}$ (= connected component of $L(\;^LG)_\phi.$)\endproclaim

\demo{Proof} We begin with uniqueness (iii): Since $\phi'= x\phi x^{-1}$ we conclude that

\noindent $ (L(\;^LG)_{\phi'},\;P(\;^LG)_{\phi'})=x(L(\;^LG)_{\phi},\;P(\;^LG)_{\phi})x^{-1},$ and this implies

\noindent $(L(\;^LG)_{\phi'},\;P(\;^LG)_{\phi'})=(L(\;^LG)_{\phi},\;P(\;^LG)_{\phi})$ because all groups are standard. Comparing the connected components we see that $x$ normalizes $P(\wh G)_\phi$ hence $x\in P(\wh G)_\phi.$ But $x$ also normalizes the Levi group $L(\wh G)_\phi$ and therefore
$$   x\in L(\wh G)_\phi,\tag{1}$$
because the unipotent radical $R_u$ of $P(\wh G)_\phi$ acts simply transitive on the Levi complements in $P(\wh G)_\phi = L(\wh G)_\phi\ltimes R_u.$ Moreover (1) implies
$$  z(\phi')= xz(\phi)x^{-1} =z(\phi),$$
because $z(\phi)$ is in the connected center of $L(\;^LG)_\phi.$

{\bf (i)/(ii)} We choose any $\phi\in [\phi]$ and consider the corresponding pair $(z(\phi),\; L(\,^LG)_\phi)$ as in {\bf 5.5}. We have $Im(\phi)\subseteq L(\,^LG)_\phi$ and so, by definition {\bf 3.1}, the Levi group $L(\,^LG)_\phi$ must be relevant which means it is $\wh G$-conjugate to a relevant standard group (see {\bf 2.5.5}). Since the assignments are compatible with $\wh G$-conjugation (see {\bf 5.5}) we may choose $\phi_*\in [\phi]$ such that $L(\;^LG)_{\phi_*}$ is a relevant standard group which according to {\bf 2.5.4}  means that

$$   L(\;^LG)_{\phi_*} =\;^LM_*,$$
where $M_*\supseteq M_0$ is a standard Levi subgroup in $G=\bG(F),$ and according to {\bf 4.3(3)} we have:
$$   z({\phi_*}) \in  Z(\;^LM_*)^0_{hyp} \cong  \aa^*_{M_*}.$$
This determines $\nu({\phi_*})\in \aa_{M_*}^*.$ But:

\noindent {\bf The euclidean space $\aa^*_{M_*}$ sees all $F$-Levi subgroups $M'\supset M_*$ as hyperplanes. (see 1.2)} And by the condition {\bf (P1)} the element $\nu({\phi_*})$ cannot lie on such a hyperplane, because otherwise $\;^LM_*$ would not be maximal with respect to $z({\phi_*}).$ (We use here that all identifications {\bf 4.3(3)} are compatible with each other in a natural way, because they are all subordinate to {\bf 4.3(4)}.) Therefore $\nu({\phi_*})$ determines a chamber in $\aa^*_{M_*}$ which must correspond to a semi-standard parabolic group $P_*=\bP(F)\subseteq G$ with Levi group $M_*.$ In particular we have $M_*=M_{P_*}$ as in {\bf 1.1}.

\noindent Now proceeding with $(P_*,M_*),$ we find $w\in \;_FW,$ in the relative Weyl group, such that
$$  (P,M) = w(P_*,M_*)w^{-1}$$
is a pair consisting of standard groups in $G.$ We note that $P$ is uniquely determined by $P_*,$ and $M_*=M_{P_*}$ implies that $M=M_{P}$ is again standard, but we may have $M\ne M_*$  because different standard Levi groups can be conjugate. Now we are going to use the $W$-$\wh W$ equivariant correspondences (5), (5a) of {\bf 2.2}.
As we have seen in {\bf 2.5.4} they induce $\;_FW$-$\;_F\wh W$ equivariant correspondences between semi-standard groups in $\bG$ which are defined over $F$, and relevant semi-standard subgroups of $\;^LG.$

\noindent Using these correspondences we can shift
$$  (P_*,M_*)\mapsto w(P_*,M_*)w^{-1} = (P, M),\quad\text{with}\quad w\in \;_FW,$$
into
$$  (\wh P_*,\wh M_*)\mapsto \wh w(\wh P_*,\wh M_*) \wh w^{-1} = (\wh P, \wh M),$$
where $w\leftrightarrow \wh w\in\;_F\wh W.$ The groups
$$   \;^LP=\wh P\rtimes W_F,\qquad \;^LM=\wh M\rtimes W_F $$
are then standard groups in $\;^LG,$ which are again relevant. Moreover $\wh w\in\;_F\wh W$ induces
$$  Z(\;^LM_*)^0\mapsto \wh w Z(\;^LM_*)^0\wh w^{-1} = Z(\;^LM)^0,\qquad \phi:= \wh n_{\wh w}\phi_*\wh n_{\wh w}^{-1},$$
where $\wh n_{\wh w}\in N_{\wh G}(\wh T)$ represents $\wh w\in N_{\wh W}(\wh W_0)^\Gamma \subseteq \wh W.$
Then $z(\phi_*)\mapsto \wh w z(\phi_*)\wh w^{-1} = z(\phi),$
and the triple $(\phi,\;^LM,\;^LP)$ is what we are looking for.
\qed\enddemo

\proclaim{6.2 Corollary} Let $[\phi]_G\in \Phi(G)$ be an L-parameter and choose a triple $(\phi,M,P)$ where $\phi\in [\phi]_G,$ and $P,$ $M=M_P$ are standard F-groups in $G=\bG(F)$ such that

$\;^LM= L(\;^LG)_\phi,$  and

$ z(\phi)\in Z(\;^LM)^0_{hyp} \leftrightarrow \nu(\phi)\in\aa_{M}^* $ is regular with respect to $P.$

\noindent Then the triple $(P, [\tphi]_M, \nu)$  where $\tphi = \phi_{z^{-1}(\phi)}$ and $\nu=\nu(\phi)$ is uniquely determined and is the L-parameter standard-triple which is assigned to $[\phi]_G.$\endproclaim

Because this map from L-parameters to L-parameter standard-triples is obviously the inverse to the map we have defined in {\bf 4.6}, the proof has been finished now.

\bigskip

{\bf 7. A remark on refined L-parameters}

Modulo a certain ambiguity to express irreducible representations by standard triples and L-parameters by L-parameter standard triples, (see 1.4 and 2.8 resp.), we fix identifications
$$   (P,\sigma_M,\nu)=\pi\in Irr(G),$$
$$   (P,[\tphi]_M,\nu)=[\phi]\in \Phi(G).$$
Now, if we assume that for all standard Levi subgroups $M$ of $G$ (including $M=G$), a tempered reciprocity map
$$   \sigma_M\in Irr(M)_{temp} \mapsto [\tphi](\sigma_M)\in \Phi(M)_{temp}\tag{1}$$
has been given, then we can extend this in a natural way to a reciprocity map for all representations:
$$   (P,\sigma_M,\nu)\in Irr(G)\quad\mapsto \quad (P,\;[\tphi](\sigma_M),\;\nu)\in \Phi(G),\tag{2}$$
where the parameters $(P,\nu)$ are left fixed!
Moreover if we fix an L-parameter $[\phi]=(P,[\tphi]_M,\nu)\in \Phi(G),$ then we have a natural identification
$$   \Pi_{[\phi]} \leftrightarrow \Pi_{[\tphi]_M},\tag{3}$$
between the fiber (or L-packet) of all representations $\pi\in Irr(G)$ which are mapped onto $[\phi]$ and the tempered L-packet of all representations $\sigma_M\in Irr(M)_{temp}$ which are mapped onto $[\tphi]_M\in\Phi(M)_{temp}.$ Obviously the bijection {\bf (3)} from right to left is obtained by associating $(P,\sigma_M,\nu)$ to the tempered representation $\sigma_M\in \Pi_{[\tphi]_M}.$ In particular the data $(P,\nu)$ are fixed for all $\pi\in\Pi_{[\phi]},$ because they are determined already by $[\phi],$ and therefore we cannot have tempered representations in $\Pi_{[\phi]}$ if the L-parameter $[\phi]$ is not tempered. (The possibility of such an identification (3) has been predicted in {\bf [A],p.201}.)

\noindent As to {\bf (3)} we mention that the members of an L-packet $\Pi_\phi$ are expected to be in correspondence with (certain) irreducible representations of the finite group
$$\overline S_\phi = C_{\wh G}(Im\;\phi)/\{C_{\wh G}(Im\;\phi)^0\cdot Z(\wh G)^\Gamma\},$$
which, due to the identity $Z(\wh G)^\Gamma = C_{\wh G}(\;^LG)$ can be represented by the exact sequence
$$  C_{\wh G}(\;^LG)/C_{\wh G}(\;^LG)^0\rightarrow C_{\wh G}(Im\;\phi)/C_{\wh G}(Im\;\phi)^0\rightarrow \overline S_\phi\rightarrow 1,\tag{4}$$
where the middle term we will call $S_\phi.$ Now concerning {\bf (3)} we have to compare the groups $\overline S_\phi(G)$ and $\overline S_{\tphi}(M).$ Since we have now
$$  Im(\phi) \subset \;^LM \subseteq\;^LG,$$
we may form two component groups $S_\phi(M)$ and $S_\phi(G),$ and obviously there is a
natural map $$  S_\phi(M)\rightarrow S_\phi(G).\tag{5}$$

\noindent But according to {\bf 5.6} the Levi group $\;^LM=L(\,^LG)_\phi$ is maximal such that {\bf 5.(P1)} holds. And via {\bf 5.(P1)*} this is equivalent to  $\;^LM =L(\;^LG)_\phi= C_{\;^LG}(\wh S),$ where $\wh S$ is minimal such that $z(\phi)\in \wh S\subset Z(C_{\wh G}(Im(\phi))^0$. So we see that
$C_{\wh G}(Im(\phi))\subset C_{\,^LG}(\wh S)=\,^LM$ which implies $C_{\wh G}(Im(\phi))=C_{\wh M}(Im(\phi)),$
hence {\bf (5)} is an equality. But actually we have
\proclaim{7.1 Proposition} If $[\phi]\in \Phi(G)$ is given by the L-parameter standard triple

\noindent $(P,[\tphi],\nu)$ then in accordance with {\bf (3)} we obtain an equality
$$  \overline S_\phi(G) = \overline S_{\tphi}(M_P).$$\endproclaim
\demo{Proof} We have $\overline S_{\tphi}(M)= \overline S_\phi(M)$ because the images of $\tphi$ and of $\phi$ in $\,^LM$ differ only by a central twist. So equivalently we have to see that $S_\phi(M)=S_\phi(G)$ induces $\overline S_\phi(M)=\overline S_\phi(G)$. For this we use Lemma 1.1 of {\bf [A 99]}:
$$  C_{\wh M}(\,^LM)= C_{\wh G}(\,^LG)\cdot C_{\wh M}(\,^LM)^0,$$
which is true for all semistandard Levi subgroups $\,^LM\subset \,^LG.$
Using this identity we can compare the $M$-version and the $G$-version of {\bf (4)} as follows:
$$  C_{\wh M}(\,^LM)\cdot C_{\wh M}(Im\,\phi)^0 =C_{\wh G}(\,^LG)\cdot C_{\wh M}(\,^LM)^0\cdot C_{\wh M}(Im\,\phi)^0,$$
and $C_{\wh M}(\,^LM)^0\subseteq C_{\wh M}(Im\,\phi)^0$ implies now that:
$$  C_{\wh M}(\,^LM)\cdot C_{\wh M}(Im\,\phi)^0=C_{\wh G}(\,^LG)\cdot C_{\wh M}(Im\,\phi)^0=C_{\wh G}(\,^LG)\cdot C_{\wh G}(Im\,\phi)^0$$
hence the result.\qed\enddemo

As a consequence we have got now identifications
$$  Irr(\overline S_{\tphi}(M_P)) = Irr(\overline S_\phi(G))\qquad CL(\overline S_{\tphi}(M_P)) = CL(\overline S_\phi(G)),\tag{6}$$
between irreducible representations and conjugacy classes resp. if $[\phi]$ is given by the L-parameter standard triple $(P,[\tphi],\nu).$

\noindent According to {\bf [A]} it is widely believed that the members of the tempered L-packet $\Pi_{[\tphi]_M}$ are in 1-1 correspondence with the irreducible representations of $\overline S_{\tphi}(M)$ and with the elements of $CL(\overline S_{\tphi}(M))$ resp. if the group $G=\bG(F)$ is quasisplit. Moreover in loc.cit. Arthur conjectures that in general the members of $\Pi_{[\tphi]_M}$ should be in 1-1 correspondence with a certain subset of
$CL(\overline S_{\tphi}(M)).$ Using the Langlands classification of L-parameters we may transport this from tempered to nontempered L-packets.

\bigskip

{\bf 8. Examples}

\noindent {\bf Example 1:} We consider $\bGL_n|F$, and as a Borel subgroup and a maximal torus we fix $\bB_n\supset \bT_n$ the upper triangular and the diagonal matrices resp. As a basis $e_1,...,e_n\in X^*(\bT_n)\cong\Z^n$ we fix the coordinate maps
$$  t=\text{diag}(t_1,...,t_n)\in\bT_n \mapsto  t^{e_i}:= t_i \in \bGL_1.\tag{1}$$
The Weyl group $W=S_n=W(\bT_n,\bGL_n)$ acts on $X^*(\bT_n)$ by permuting the rational characters $e_i,$ and therefore the pairing
$$  \langle e_i,e_j \rangle =\delta_{i,j}$$
induces a Weyl-group-invariant euclidean structure an $\aa_{\bT_n}^*= \R\otimes X^*(\bT_n).$ The simple roots $\Delta(\bT_n,\bB_n)\subset X^*(\bT_n)$ are:
$$  \Delta(\bT_n,\bB_n)=\{e_i-e_{i+1}\;|\; i=1,...,n-1\}, \qquad t\mapsto t^{e_i-e_{i+1}}= t_i/t_{i+1}.$$
For the adjoint action of $\bT_n$  on $\text{Lie}(\bGL_n)$ the position $(i,i+1)$ yields the eigenspace of the rational character $e_i-e_{i+1}.$ So we identify the simple roots with the corresponding positions. At the same time the transposition $(i,i+1)\in W=S_n$ is the reflection corresponding to $e_i-e_{i+1}.$

\noindent If $I\subseteq \Delta(\bT_n,\bB_n)$ is a subset then the standard Levi group $\bM_I=\bM(W_I)\supseteq\bT_n$ is the minimal bloc-diagonal subgroup of $\bGL_n$ which includes all the positions $(i,i+1)\in I.$ With respect to $\bM_I$ we obtain
$$  I= \Delta(\bT_n,\bB_n\cap\bM_I),\qquad W_I = W(\bT_n,\bM_I),$$
and $\bP_I=\bM_I\bB_n= \bB_n W_I\bB_n$ is the corresponding standard parabolic subgroup.

\noindent If $I'=\Delta(\bT_n,\bB_n)-I$ is the complement, and if $I'$ corresponds to the positions
$$I'=\{(i_1,i_1+1),...,(i_s,i_s+1)\},\qquad 0< i_1< i_2 <...< i_s <n,$$
then we obtain
$$ \bM_I = \prod_{\k =1}^{s+1} \bGL_{m_\k},\quad m_\k=m_\k(I):=i_\k - i_{\k-1},\quad \sum_{\k=1}^{s+1} m_\k =n,\tag{2}$$
with the notation $i_0=0,\;i_{s+1}=n.$

\noindent Let $b_1,...b_{s+1}\in X^*(\bM_I)$ be the basis where
$b_\k :\bM_I\rightarrow \bGL_1$ is the determinant on $\bGL_{m_\k}$ and is trivial on the other factors. Then we have
$$ res_{\bT_n}(b_\k) = e_{i_{\k-1} +1}+\cdots +e_{i_\k}\in X^*(\bT_n), \quad\text{for all}\;\k=1,...,s+1.$$
In the following we will identify
$$ b_\k = e_{i_{\k-1} +1}+\cdots +e_{i_\k}\in X^*(\bM_I)\hookrightarrow X^*(\bT_n),\tag{3}$$
and we will consider $\{b_1,...,b_{s+1}\}$ as a basis of the euclidean space $\aa_{\bM_I}^*=\R\otimes X^*(\bM_I)$ of dimension $s+1=n-\# I.$

Let $\bA_I$ be the central torus of $\bM_I.$ Since $\bA_I=\prod_{\k=1}^{s+1} \overline F^\times \cdot I_{m_\k}\subseteq \bT_n,$ we obtain the surjection
$$ e_i\in X^*(\bT_n)\mapsto \o e_\k\in X^*(\bA_I),$$
where  $\o e_\k := res(e_i)$ for any $i$ such that $i_{\k-1}<i\le i_\k.$ The combined map is: $$X^*(\bM_I)\rightarrow X^*(\bA_I),\qquad b_\k \mapsto e_{i_{\k-1} +1}+\cdots +e_{i_\k}
\mapsto m_\k\cdot\o e_\k.$$
\bigskip

\noindent {\bf Now we identify} $\aa_{\bM_I}^*=\R\otimes X^*(\bM_I)=\R\otimes X^*(\bA_I),$ hence
$$  \o e_\k=\frac{1}{m_\k} b_\k\quad \in\aa_{\bM_I}^*\subseteq \aa_{\bT_n}^*,$$
$$  <\o e_\k,\o e_\k>\;=\; <\frac{1}{m_\k}b_\k,\frac{1}{m_\k}b_\k>\; =\; \frac{1}{m_\k},\qquad
<b_\k,\o e_\k> =1. $$

\noindent The elements from $I'$ (which are the simple roots outside of $\bM_I$) give rise to roots $e_{i_\k}-e_{i_\k +1}\in X^*(\bT_n),$ such that
$$ res(e_{i_\k}-e_{i_\k +1}) =\o e_\k -\o e_{\k+1}\in X^*(\bA_I),$$
form the root basis $\Delta(\bP_I)$ with respect to the adjoint action of $\bA_I$ on the unipotent radical $U(\bP_I).$ Now we see from {\bf 1.3 c)} that:
\proclaim{8.1 Lemma} The regularity of $$\nu= \sum_{\k=1}^{s+1}\beta_\k b_\k\quad  \in\aa_{\bM_I}^*= \R\otimes X^*(\bA_I)$$
with respect to our parabolic group $\bP_I$ means that
$$  <\nu,\o e_\kappa -\o e_{\kappa+1}>= \beta_\k -\beta_{\k+1}\quad > 0\quad \text{for all}\;\k=1,...,s,$$
or equivalently $\beta_1 >\cdots >\beta_{s+1}$ for the coefficients of $\nu$ with respect to the determinant characters $b_1,...,b_{s+1}.$\endproclaim

Using this Lemma the Langlands classification {\bf 1.4} for $G=GL_n(F)$ has the following well known reformulation (see for instance {\bf [Ku], Theorem 2.2.2.}):
\proclaim{8.2 Langlands classification for $GL_n(F)$} The irreducible representations $\pi$ of $GL_n(F)$ are in bijection with all possible tuples
$$\TTT = \left( m_1,...,m_{s+1};\;\tau_1,...,\tau_{s+1};\;\beta_1>\cdots >\beta_{s+1} \right)$$
where $m_i\ge 1$ are natural numbers with sum $=n,$ where $ \tau_i$ is an irreducible tempered representation of $GL_{m_i}(F)$ and where the $\beta_i$ are real numbers. The representation $\pi=\pi(\TTT)$ is obtained as
$$  \pi(\TTT)= j(i_{GL_n(F),P(\TTT)}\left(\tau_1\cdot|det|_F^{\beta_1}\otimes\cdots\otimes \tau_{s+1}\cdot|det|_F^{\beta_{s+1}}\right))$$
where $P(\TTT)$ is the standard parabolic group with Levi subgroup $M(\TTT)=GL_{m_1}(F)\times\cdots\times GL_{m_{s+1}}(F)$, and where $\pi(\TTT)$ is the unique irreducible quotient of the normalized parabolic induction.\endproclaim
For the proof one only needs to observe that $M(\TTT)= M_I$ in our previous notation, and the tempered representations of $M_I$ come as outer tensor products of tempered representations for the different $GL_{m_i}(F)$. And under
$\aa_{\bM_I}^*\rightarrow X_{ur}(M_I)_+$ the regular element $\nu=\beta_1b_1+\cdots\beta_{s+1}b_{s+1}$ is mapped to the unramified character
$$  \chi_\nu = (|det|_F^{\beta_1}\otimes\cdots\otimes |det|_F^{\beta_{s+1}})\in X_{ur}(M_I)_+,$$
as we see from the remarks following {\bf 1.(1)}. With this in mind, 8.2 comes as a direct application of 1.4.

Before we can come to the Langlands classification of L-parameters, we need to consider also the dual group
$$ (\wh G= GL_n(\C),\;\wh B= B_n(\C),\;\wh T= T_n(\C)).$$
The identification $X_*(\wh T) =X^*(\bT_n)$ is obtained via $   \wh e_i = e_i,$ where
$$\wh e_i:\C^\times\rightarrow \wh T,\qquad \lambda^{\wh e_i}=\text{diag}(1,...,\lambda,1,...)\tag{4}$$
maps $\lambda$ to the diagonal position $(i,i).$ The coroots $\wh e_i - \wh e_{i+1}\in \Delta_*(\wh T,\wh B)$ identify with the roots $e_i-e_{i+1}\in\Delta(\bT_n,\bB_n).$

\noindent From $\wh T\cong \C^\times\otimes X_*(\wh T)$ and $X_*(\wh T)=X^*(\bT_n),$ we recover the hyperbolic elements as
$$  \wh T_{hyp}\cong(\R_+)^\times\otimes X_*(\wh T)\cong \R\otimes X^*(\bT_n)=\aa_{\bT_n}^*\tag{5}$$
where from right to left we use the exponential map $s\in\R \mapsto q^s,$ hence

\noindent $q^{s\otimes \wh e_i}=(q^s)^{\wh e_i}\in\wh T_{hyp}$ corresponds to $s\otimes e_i\in \aa_{\bT_n}^*.$

Now we go back to our standard Levi subgroup $\bM_I\subset \bGL_n$ and turn to the dual group $\wh{\bM_I}=\wh M_{\wh I}\subset \wh G.$ From {\bf (2)} it is obvious that
$$  \wh M_{\wh I}  = \prod_{\k =1}^{s+1} GL_{m_\k}(\C) \subset \wh G=GL_n(\C).$$
And from {\bf (3), (4)} we see that $b_\k= e_{i_{\k-1} +1}+\cdots +e_{i_\k}\in X^*(\bM_I)\subset X^*(\bT_n)$ must correspond to
$\wh b_k=\wh e_{i_{\k-1} +1}+\cdots +\wh e_{i_\k}\in X_*(Z(\wh M_{\wh I}))$ such that
$$  \lambda\in\C^\times \mapsto \lambda^{\wh b_\k} = \text{diag}\left(I_{m_1},\;...,\;\lambda I_{m_\k},\; I_{m_{\k+1}},\;...\right) \in Z(\wh M_{\wh I}).$$
{\bf In particular} $b=det\in X^*(G)$ corresponds to $\wh b=\{\lambda\mapsto \lambda I_n\}\in X_*(Z(\wh G)).$
Therefore:
\proclaim{8.3 Lemma} The bijection {\bf (5)} restricts to a bijection $z(\nu)\leftrightarrow \nu:$
$$  Z(\wh M_{\wh I})_{hyp}\ni \text{diag}\left(q^{\beta_1}I_{m_1},\;...,\;q^{\beta_{s+1}} I_{m_{s+1}}\right) \longleftrightarrow \sum_\k\beta_\k b_\k \in\aa_{\bM_I}^*\tag{6}$$
and the regularity of $\nu=\sum_\k\beta_\k b_\k$ is equivalent to $q^{\beta_1}>\cdots >q^{\beta_{s+1}} >0.$
\endproclaim

Now we are prepared to give the Langlands classification of L-parameters for
$G= GL_n(F).$ This implies $\,^LG= GL_n(\C)\times W_F$ a direct product because the $\Gamma$-action on $\wh G= GL_n(\C)$ is trivial. Therefore in {\bf 3.1 (i)} we may replace $\phi=(\phi_1,\phi_2)\in \Phi(G)$ by $(\phi_1,\vf_2)$ where $\vf_2:W_F\rightarrow GL_n(\C)$ is now a continuous and Frobenius-semisimple homomorphism, and we obtain
$$  \Phi(G) = \text{Hom}_0(SL_2(\C)\times W_F,\; GL_n(\C))/\sim.$$
The L-parameter $\phi$ is tempered if and only if $Im(\vf_2)$ is relatively compact. Again we consider standard Levi subgroups $M_I=\bM_I(F)\subseteq G.$ Then we have $M_I=\prod_\k GL_{m_k}(F),$ $\wh{M_I}=\prod_\k GL_{m_k}(\C)$ and therefore
$$  \Phi(M_I) = \text{Hom}_0(SL_2(\C)\times W_F,\; \wh{M_I})/\sim  \quad =\prod_{\k=1}^{s+1} \Phi(GL_{m_\k}(F)).$$
Moreover it is obvious that $\phi=(\phi^{(1)},...,\phi^{(s+1)})\in \Phi(M_I)$ is tempered if and only if all components $\phi^{(i)}$ are tempered L-parameters. Finally with $s\in\C$  we associate the unramified character
$$  \omega_s:=\{w\in W_F \mapsto \pa w\pa^s\}$$
(see {\bf 4.7} and {\bf [T],(2.2)} resp.). Then from {\bf 8.3} and {\bf 4.6} we obtain the following
\proclaim{8.4 Langlands classification of L-parameters for $GL_n(F)$} The L-parameters $[\phi]\in\Phi(GL_n(F))$ are in bijection with all possible tuples
$$ \TTT = (m_1,...,m_{s+1};\;[\tphi^{(1)}],...,[\tphi^{(s+1)}];\;\beta_1>\cdots >\beta_{s+1})$$
where $m_i\ge 1$ are natural numbers with sum $=n$, where $[\tphi^{(i)}]$ is a tempered L-parameter for $GL_{m_i}(F),$ and where the $\beta_i$ are real numbers. The L-parameter $[\phi]=[\phi](\TTT)$ is obtained as
$$  [\phi(\TTT)] = [\left([\tphi^{(1)}\omega_{\beta_1}],...,[\tphi^{(s+1)}\omega_{\beta_{s+1}}]\right)],$$
which means it is the uniquely determined L-parameter of $GL_n(F)$ which is generated by the given L-parameter-tuple for the standard Levi group $M(\TTT)=GL_{m_1}(F)\times\cdots\times GL_{m_{s+1}}(F).$
\endproclaim
\demo{Proof} Put $M=M(\TTT)\subseteq GL_n(F)$ and $P=P(\TTT)$ the corresponding standard parabolic group. By {\bf 8.1} we know that $\nu=\sum_\k \beta_\k b_\k\in \aa_M^*$ is regular with respect to $P.$ And $[\tphi_M]:=([\tphi^{(1)}],...,[\tphi^{(s+1)}])$ is a tempered L-parameter for $M.$ Furthermore from {\bf 8.3} we obtain
$$  z(\nu) = \text{diag}(q^{\beta_1}I_{m_1},...,q^{\beta_{s+1}}I_{m_{s+1}})\in Z(\wh M),$$
and the twist $(\tphi_M)_{z(\nu)}$ in the sense of {\bf 4.6} yields
$$ (\tphi_M)_{z(\nu)} = (\tphi^{(1)}\omega_{\beta_1},...,\tphi^{(s+1)}\omega_{\beta_{s+1}})$$
as we see from the definition of $\omega_s.$ Therefore our assertion is now a direct consequence of {\bf 4.6}.
\qed\enddemo
\bigskip

\noindent{\bf Example 2:} Now we consider inner forms of $\bG=\bGL_n$ which are the groups $\bG'|F$ such that
$$  G'=\bG'(F) = GL_m(D),$$
where $D|F$ is a fixed central division algebra of index $d=\frac{n}{m},$ in particular $m\mid n.$ Then we have
$$  X^*(G')=X^*(\bG') = \Z\cdot Nrd$$
generated by the reduced norm which comes by restricting the determinant of $\bG'\cong \bGL_n$ to the rational points $G'.$ The dual group is $\wh{G'} =GL_n(\C)$ with center $Z(\wh{G'})= GL_1\otimes I_n.$ The identification
$$  X_*(Z(\wh{G'})) = X^*(G'),\qquad Nrd^\vee \leftrightarrow Nrd $$
is defined by $Nrd^\vee(\lambda)=\lambda\cdot I_n.$  Since $(\bG',\bB',\bT')$ is an inner form of the split group $(\bGL_n(F),\bB_n,\bT_n)$ we know from {\bf 2.2.1} and the subsequent remark that the $\mu_{\bG'}(\Gamma)$-action on $X^*(\bT')$ is trivial, and therefore according to {\bf 2.3.3} the usual $\Gamma$-action on $X^*(\bM')$ will be trivial if $\bM'\supset\bT'$ is a semi-standard $F$-Levi subgroup, hence $X^*(M')=X^*(\bM').$ The minimal $F$-Levi subgroup is $\bM_0'\cong \bGL_d^{\times m}$ such that $\bM_0'(F)= GL_1(D)^{\times m}.$ The corresponding Weyl group $W_0=W(\bT',\bM_0')$ is a product of $m$ copies of the symmetric group $S_d,$ and the normalizer $N_W(W_0)$ in $W=W(\bT',\bG')\cong S_n$ are the permutations which respect the partition into subsets of $d$ elements, hence the relative Weyl group $\,_FW=N_W(W_0)/W_0\cong S_m$ comes as the group of permutations of the $m$ different subsets of order $d,$ and it acts now on $X^*(M'_0)\cong \Z\cdot e_1 +\cdots \Z\cdot e_m$ with the following rational characters $e_i:$
$$  m_0'=(d_1,...,d_m)\in D^\times \times\cdots\times D^\times \mapsto {m_0'}^{e_i}:= Nrd_{D|F}(d_i)\in F^\times,$$
where the reduced norms stand for the different copies of $D|F.$ In analogy to example 1, the definition
$$  \langle e_i, e_j\rangle= \delta_{i,j} $$
fixes a $\,_FW$-invariant euclidean structure on $\aa_{M_0'}^*=\R\otimes X^*(M_0').$

\noindent If we consider now a standard $F$-Levi subgroup
$$  M'=\bM'(F)=\prod_{\k=1}^r GL_{m_\k}(D),\qquad \sum_\k m_\k =m,\qquad \wh{M'} =\prod_{\k=1}^r GL_{m_\k d}(\C),$$
then $X^*(M')$ is the rank-$r$-lattice generated by the reduced norms $Nrd_\k$ on $GL_{m_\k}(D)$ trivially extended to the other factors, and $X_*(Z(\wh{M'}))$ is the lattice generated by the $Nrd_\k^\vee$ which are the cocharacters:
$$ Nrd_\k^\vee: \lambda\in\C^\times \mapsto \text{diag}(I_{m_1d},...,\lambda I_{m_\k d},...).$$
Using the natural embedding $X^*(M')\hookrightarrow X^*(M_0')$ we obtain the identification
$$   Nrd_\k = e_{i_{\k-1}+1}+\cdots + e_{i_\k},$$
and similar considerations as in example 1 will give:
\proclaim{8.5 Lemma} The regularity of
$$ \nu = \sum_{\k=1}^r \beta_\k Nrd_\k\in\R\otimes X^*(M'),\qquad M'=\prod_{\k=1}^r GL_{m_\k}(D),$$
with respect to the corresponding standard parabolic group $P'\supseteq P_0'$ means that the real coefficients $\beta_1 >\cdots > \beta_r$ are ordered.\endproclaim

\noindent {\bf As to the L-parameters} we may identify $\,^LG'=\,^LG= GL_n(\C)\times W_F$ and therefore as in example 1 we may switch to L-homomorphisms
$$  \phi=(\phi_1,\vf_2): SL_2(\C)\times W_F \rightarrow GL_n(\C)=\wh{G'}$$
where the values are in the connected component of $\,^LG'.$ So the L-parameters are basically the same as in the case of $GL_n(F).$ But if $d=\frac{n}{m} >1 $ then $\bG'$ is nonsplit and therefore according to {\bf 3.1} a Levi subgroup $L(\wh{G'})$ need to be {\bf relevant} if it contains $Im(\phi).$ If we assume (modulo conjugation) that the Levi subgroup $L(\wh{G'})\subseteq GL_n(\C)$ is standard, then this means that all diagonal blocks must be of size $m_i d\times m_i d$ a multiple of $d.$ Any other L-homomorphisms with image in smaller Levi-subgroups are not allowed.

\proclaim{8.6 Langlands classification of L-parameters for $GL_m(D)$} The L-parameters $[\phi]\in\Phi(GL_m(D))$ are in bijection with all possible tuples
$$ \TTT = (m_1,...,m_{s+1};\;[\tphi^{(1)}],...,[\tphi^{(s+1)}];\;\beta_1>\cdots >\beta_{s+1})$$
where $m_i\ge 1$ are natural numbers with sum $=m$, where $[\tphi^{(i)}]$ is a tempered L-parameter for $GL_{m_i}(D),$ and where the $\beta_i$ are real numbers. The L-parameter $[\phi]=[\phi](\TTT)$ is obtained as
$$  [\phi(\TTT)] = [\left([\tphi^{(1)}\omega_{\beta_1}],...,[\tphi^{(s+1)}\omega_{\beta_{s+1}}]\right)],$$
which means it is the uniquely determined L-parameter of $GL_m(D)$ which is generated by the given L-parameter-tuple for the standard Levi group

\noindent $M(\TTT)=GL_{m_1}(D)\times\cdots\times GL_{m_{s+1}}(D).$
\endproclaim
\noindent Note here that the tempered L-parameter $[\tphi^{(i)}]$ takes its values in $GL_{m_id}(\C),$ and we multiply it with the scalar function $\omega_{\beta_i}$ on $W_F.$
\bigskip

{\bf Appendix:} In both examples we know much more on the L-parameters because we may think of an L-homomorphism $\phi$ as of an appropriate $n$-dimensional representation of $SL_2(\C)\times W_F.$ Therefore it can be written as
$$   \phi = \sum_i  r_i\otimes \rho_i,$$
where $r_i,$ $\rho_i$ are both irreducible rational / Frobenius-semisimple representations of $SL_2(\C)$ and $W_F$ resp. Since up to equivalence the group $SL_2(\C)$ has precisely one irreducible representation $r_n$ of dimension $n,$ we may write
$$   \phi=\sum_i r_{n_i}\otimes \rho_i,\qquad  \sum_i n_i\text{dim}(\rho_i) = n,$$
and in case of example 2 we must have $d\mid n_i\text{dim}(\rho_i) $ for all $i.$ Furthermore every irreducible representation $\rho$ of $W_F$ can be written as
$$  \rho=\rho_0\otimes\omega_s,$$
where $\rho_0$ is of Galois type, in particular it is unitary, and $\omega_s=\{w\mapsto {\pa w\pa}^s\}$ is as before. (see for instance {\bf [T],(2.2.1)}.) Therefore we see from (2) that the hyperbolic part
$z(\phi)= \vf_2(w_1)_{hyp}$ is the diagonal matrix
$$z(\phi) = \text{diag}(...,q^{Re(s_i)},...,q^{Re(s_i)},...)$$
where $|\omega_{s_i}| = \omega_{Re(s_i)}$ corresponds to the irreducible representation $\rho_i,$ and where $n_i\cdot \text{dim}(\rho_i)$ is the number of occurences of the entry $q^{Re(s_i)}.$ A parameter $\phi$ is tempered if and only if $\vf_2$ is relatively compact, i.e. $Re(s_i)=0$ for all $i.$ Similar if we consider parameters of Levi-subgroups which come as ordered sets of parameters for the different factors $GL_{n_i}(F)$ and $GL_{m_i}(D)$ resp.

Since each $r_{n_i}\otimes \rho_i$ can be reinterpreted as a segment and each $\phi$ as a multisegment one can recover everything in terms of segments and multisegments resp. We stop here because our aim was only to demonstrate the reduction to tempered representations.

\bigskip

\noindent Allan~J.~Silberger,~Dept.~of~Mathematics,~Cleveland~State~University

\noindent Home~Address:~17734~Lomond~Blvd,~Shaker~Hts,~OH~44122,~USA

\noindent email:~a.silberger\@att.net\medskip

\noindent
Ernst-Wilhelm Zink,

\noindent
Institut f\"ur Mathematik,
Humboldt-Universit\"at zu Berlin,

\noindent
Unter den Linden 6,
D-10099 Berlin,
Germany

\noindent
e-mail: zink\@math.hu-berlin.de

\end{document}